\documentclass{article}
\PassOptionsToPackage{numbers,sort&compress}{natbib}

\usepackage[utf8]{inputenc} 
\usepackage[T1]{fontenc}    
\usepackage{hyperref}       
\usepackage{url}            
\usepackage{booktabs}       
\usepackage{amsmath, amssymb, amsthm, amsfonts}       
\usepackage{nicefrac, xfrac}       
\usepackage{microtype}      
\usepackage{xcolor}         
\usepackage{bm}
\usepackage{enumitem}
\usepackage{graphicx}
\usepackage{algorithm}
\usepackage{algorithmic}
\usepackage{natbib}[numbers,sort&compress]
\usepackage{geometry}
\usepackage{authblk,textcomp}

\usepackage[cal=euler]{mathalfa}
\usepackage{libertine}
\usepackage{mathtools}

\makeatletter
\newcommand{\otherlabel}[2]{\protected@edef\@currentlabel{#2}\label{#1}}
\makeatother
\def\mat#1{#1}
\renewcommand{\vec}[1]{\boldsymbol{#1}}

\DeclareMathOperator{\extr}{extr}
\DeclareMathOperator*{\argmin}{arg\,min}

\DeclareMathOperator{\tr}{tr}

\def\dd{\mathrm{d}}

\DeclareMathOperator{\erf}{erf}
\DeclareMathOperator{\prox}{prox}

\theoremstyle{plain}
\newtheorem{theorem}{Theorem}[section]
\newtheorem{proposition}[theorem]{Proposition}
\newtheorem{lemma}[theorem]{Lemma}
\newtheorem{corollary}[theorem]{Corollary}
\theoremstyle{definition}

\newtheorem{assumption}{Assumption}
\newtheorem{model}[theorem]{Model}
\theoremstyle{remark}
\newtheorem{remark}[theorem]{Remark}

\usepackage[export]{adjustbox}
\hypersetup{pdfauthor={IdePHICS},pdftitle={GMMUniversality},%
            colorlinks, linktocpage=true, pdfstartpage=1, pdfstartview=FitV,%
    breaklinks=true, pdfpagemode=UseNone, pageanchor=true, pdfpagemode=UseOutlines,%
    plainpages=false, bookmarksnumbered, bookmarksopen=true, bookmarksopenlevel=1,%
    hypertexnames=true, pdfhighlight=/O,%
    urlcolor=orange, linkcolor=blue, citecolor=blue
        }

\title{Are Gaussian data all you need? Extents and limits of universality \\ in high-dimensional generalized linear estimation}

\author[1]{Luca Pesce}
\author[1]{Florent Krzakala}
\author[2]{Bruno Loureiro}
\author[1]{Ludovic Stephan}
\affil[1]{\small Ecole Polytechnique F\'{e}d\'{e}rale de Lausanne (EPFL). 
Information, Learning and Physics (IdePHICS) lab. \newline CH-1015 Lausanne, Switzerland.}
\affil[2]{\small Département d’Informatique, École Normale Supérieure - PSL \& CNRS, 45 rue d’Ulm, \newline F-75230 Paris cedex 05, France.}
\affil[ ]

\date{}

\begin{document}

\maketitle

\begin{abstract}
In this manuscript we consider the problem of generalized linear estimation on Gaussian mixture data with labels given by a single-index model. Our first result is a sharp asymptotic expression for the test and training errors in the high-dimensional regime. Motivated by the recent stream of results on the Gaussian universality of the test and training errors in generalized linear estimation, we ask ourselves the question: "\emph{when is a single Gaussian enough to characterize the error?}". Our formula allow us to give sharp answers to this question, both in the positive and negative directions. More precisely, we show that the sufficient conditions for Gaussian universality (or lack of thereof) crucially depend on the alignment between the target weights and the means and covariances of the mixture clusters, which we precisely quantify. In the particular case of least-squares interpolation, we prove a strong universality property of the training error, and show it follows a simple, closed-form expression. Finally, we apply our results to real datasets, clarifying some recent discussion in the literature about Gaussian universality of the errors in this context. 
\end{abstract}

\section{Introduction}
\label{sec:main:intro}

It is commonsense in machine learning that structure in the data is an important ingredient for successful learning. Quantifying this statement, and in particular how structure in the features impact the training and generalization errors the most, is an important endeavor in the broad program of "seeing through" the modern machine learning black box. On the theoretical side, there has been some important recent progress in this direction in the context of generalized linear estimation. For instance, a recent line of work on linear regression trained on Gaussian data has shown that good generalization can arise even in the "overparametrized regime" where the training error is exactly zero \cite{Bartlett_2020, hastie_2022_surprises, wu20}. This \emph{benign overfitting} property crucially depends on the covariance structure, occurring when the signal components of the target align with a lower-dimensional of the data, leaving space for the noise to spread along the higher-dimensional orthogonal subspace \cite{Bartlett_2020}. Analogous conclusions hold, under similar conditions, to generalized linear tasks \cite{JMLR:v22:20-603, Wang2021, shamir22a}. Indeed, this is only one example of many surprising insights learned from the study of generalized linear models on Gaussian data over the past few years \cite{gardner1988optimal, krogh1992simple, donoho2009observed, candes2020phase}. Despite the seemingly constraining assumption on the distribution of the features, a recent line of work provides strong evidence for the \emph{Gaussian universality} of the training and generalization errors in generalized estimation in different settings. These includes rigorous results for non-Gaussian designs \cite{8006947, Panahi2017}, random feature maps \cite{mei19, gerace20a, Goldt2021TheGE, Hu2020UniversalityLF}, neural tangent features \cite{Montanari2022UniversalityOE}, Gaussian mixtures with random labels \cite{randlab}, as well as extensive numerical evidence for other feature maps \cite{Goldt2021TheGE, gcm_paper} and even real datasets \cite{gcm_paper, Jacot2020, bordelon20}. These works beg the question \textit{"when are Gaussian features a good model for learning?"}. Our aim is to give precise answers to this question in the context of generalized linear estimation on a popular model for multi-modal data, known to be able to approximate any distribution: the Gaussian mixture model.

\subsection{Main results}
Our \textbf{main contributions} in this work are as follows:
\begin{itemize}[wide=1pt, topsep=0pt]
        \item \textbf{Exact asymptotics of GLMs:} We provide the exact asymptotic limit for the training and test errors of a generalized linear model with convex  loss in high dimensions, when the data is drawn from a Gaussian mixture model with a single-index target.
These asymptotics are based on the so-called \emph{replica method} from statistical physics, and follow the same line as \cite{gaussian_mix_group}, which considered instead the task of learning the mixture labels.

\item \textbf{Universality of training and test errors:} We provide a set of sufficient conditions on the target weights $\theta_0$ such that both the asymptotic training and test errors for a Gaussian mixture model are independent from the cluster means. In particular, these conditions are satisfied by a target whose direction is uniform on $\mathbb{S}^{d-1}$. In the  case of ridge regression, we show an even stronger result: namely, the training loss is also independent from the cluster covariances, and reduces to that of a single Gaussian with identity covariance.

\item \textbf{The importance of Homoscedasticity}: In the particular case of a homoscedastic Gaussian mixture (a mixture of Gaussians that share the same covariance matrix), we further demonstrate universality results that can actually be observed on real data after a random feature map (see e.g. Fig.~\ref{fig:fig1} and \ref{fig:fig2}). We also unveil the universal behavior of the linear separability transition, a phenomenon studied in detail for pure Gaussian data in \cite{candes2020phase} and that appears to be universal for a homoscedastic mixture.

\item \textbf{Breaking universality:} In contrast to the results of the previous paragraph,
we show that there are two ways to break Gaussian universality. First, strongly heteroscedasticity can break the universal behavior. Second, in the homoscedastic case we show that the correlation between the data and the task matters: even a small correlation between the target weights and the cluster means suffice to break universality, in the sense that the asymptotic errors of a model trained in a Gaussian mixture differs from the one of a model trained on Gaussian data. Rather than the structure of the data itself, what appears to matter is thus the correlation between this structure and the task to be learned. 
\end{itemize}

\subsection{Related works}

\paragraph{Exact asymptotics:} An appealing feature of Gaussian data is that the asymptotic performance of different models can be sharply characterized in the proportional high-dimensional limit where $n,d\to\infty$ at fixed \emph{sample complexity} $\alpha \coloneqq \sfrac{n}{d}$. This is particularly the case for ridge regression, because of the close connection to a random matrix theory problem; see e.g. \cite{dobriban_2018_high, hastie_2022_surprises, wu20}. Beyond the quadratic case, there exists many asymptotic rigorous studies, for instance \cite{thrampoulidis_2018_precise} studied M-estimators, \cite{celentano20}  the performance of the LASSO estimator, while \cite{gcm_paper} provided a general result for convex losses and penalties with arbitrary covariances.

In the case of Gaussian mixtures, most of the effort has been geared towards classification, i.e. recovering the cluster label instead of teacher-generated ones. For binary classification, examples include \cite{mai.liao.ea_2019_large, mignacco20, deng.kammoun.ea_2022_model}, the latter of which also shows an equivalence between classification and a single-index model. In the multi-class setting, \cite{thrampoulidis.oymak.ea_2020_theoretical} studied the performance of ridge regression classifiers; the most general result in this line is \cite{gaussian_mix_group}, which considers any convex (not necessarily separable) loss.

\paragraph{Gaussian universality:} Remarkably, this model is also able to capture the errors of particular classes of non-Gaussian features. This \emph{Gaussian universality property} (GEP) \cite{Goldt2020ModellingTI} was proven to hold for generalized linear estimation with random features \cite{mei19, gerace20a, Goldt2021TheGE, Hu2020UniversalityLF, Schroder2023}, neural tangent features \cite{Montanari2022UniversalityOE} and kernel features \cite{karoui2008spectrum, bordelon20, Mei2021GeneralizationEO, Cui_2022, Cui2022}. \cite{Ba2022} showed that Gaussian universality is preserved on the random features model when the weights are trained at order one steps, but break if an extensive number of steps are taken. Beyond the realm of theorems, \cite{gcm_paper} provided numerical evidence of Gaussian universality for a broader class of realistic features from trained neural networks. On a close line,  \cite{randlab} studied the Gaussian universality for pure random binary labels, which was an important source of inspiration for the present work, while \cite{bordelon20, spigler2020asymptotic,Jacot2020, gcm_paper, Cui_2022,Wei2022} has numerically shown that for ridge regression in particular, the Gaussian formula captured the learning curves of some simple real datasets.

\section{Setting \& motivation}
\label{sec:main:model}
Let $(\vec{x}^{\nu}, y^{\nu})\in\mathbb{R}^{d}\times \mathcal{Y}$ denote $\nu\!=\!1,\cdots, n$ pairs of independently sampled training points. We shall be interested in studying the properties of generalized linear estimation $\hat{y}(\vec{x}) = \hat{f}(\vec{\theta}^{\top}\vec{x})$ with weights $\vec{\theta}\in\mathbb{R}^{d}$ learned from the training data by minimizing the following empirical risk:
\begin{align}
    \hat{\mathcal{R}}^{\lambda}_{n}(\vec{\theta}) = \frac{1}{n}\sum\limits_{\nu=1}^{n}\ell\left(y^{\nu}, \frac{\vec{\theta}^{\top}\vec{x}}{\sqrt{d}}\right) + \lambda r(\vec{\theta})
    \label{eq:main:erm}
\end{align}
\noindent where $\ell:\mathcal{Y}\times \mathbb{R}\to \mathbb{R}_{+}$ is a convex loss function and $r:\mathbb{R}^d\to \mathbb{R}$ is a convex penalty. For example, this includes the particular case of ridge regression  where $\ell(y,\hat{y}) = (y-\hat{y})^2$ and $r(\Vec{\theta}) = ||\Vec{\theta}||_2^2$. The key quantities of interested in the following will be the training and generalization error, defined as:
\begin{align}
    \label{eq:main:erm_train_err}
    \varepsilon_{\text{tr}} (\hat{\Vec{\theta}}) &= \frac{1}{n}\sum\limits_{\nu=1}^{n}\ell \left(y^{\nu}, \frac{\hat{\vec{\theta}}^{\top}\vec{x}^{\nu}}{\sqrt{d}}\right) \\ 
    \label{eq:main:erm_gen_err}
    \varepsilon_{\text{gen}}(\hat{\Vec{\theta}}) &= \mathbb{E}\left[
    g\left( y_{\text{new}},\hat{f}\left(\frac{\hat{\Vec{\theta}}^\top \Vec{x}_{\text{new}}}{\sqrt{d}}\right)\right)
    \right]
\end{align}
\noindent where $g:\mathcal{Y}\times \mathbb{R}\to\mathbb{R}_{+}$ is a performance metric of the choice of the statistician, not necessarily equal to the loss function $\ell(y,\hat{y})$. For example, in the case of binary classification with $\mathcal{Y}=\{-1,+1\}$, we can take $\ell(y,\hat{y}) = \log(1+e^{-y\hat{y}})$ to be the logistic loss, while taking $g(y,\hat{y}) = \mathbb{I}\left[y \neq  \hat{y}\right]$ to be the classification error. Note that in eq.~\eqref{eq:main:erm_gen_err} the expectation is taken over a new data pair $(\vec{x}_{\text{new}}, y_{\text{new}})$ which we assume is independently drawn from the same distribution of the training data.

In particular, we will be interested in characterizing these errors under the assumption that the labels have been generated by the following target distribution:
\begin{align}
\label{eq:label_generation}
	y^{\nu} \sim P_0\left(\cdot \,\big\vert\, \frac{\vec{\theta_0}^\top \vec{x}^{\nu}}{\sqrt{d}}\right)
\end{align}
\noindent for some fixed vector $\vec{\theta}_{0}\in\mathbb{R}^{d}$ and distribution $P_0$. A common choice for $P_0$ is $y^\nu = f_0\left(\frac{\vec{\theta_0}^\top \vec{x}^{\nu}}{\sqrt{d}} + \xi \right)$ for additive Gaussian noise $\xi\sim\mathcal{N}(0,\Delta)$. This setting is sometimes refereed to as a \emph{teacher-student setting}. We will sometimes adopt this convenient terminology, refereeing to the target distribution $P_0$ as the \emph{teacher} and $\vec{\theta}_{0}\in\mathbb{R}^{d}$ as the \emph{teacher weights}. Similarly, we will sometimes refer to the model $\hat{y}$ as the \emph{student} and $\hat{\vec{\theta}}\in\mathbb{R}^{d}$ the \emph{student weights}. As previously mentioned, we shall be considering both regression $\mathcal{Y}=\mathbb{R}$ and binary classification $\mathcal{Y} = \{-1,+1\}$.

This model has been the subject of a plethora of works in the high-dimensional statistics literature over the past few years, in particular under the \emph{Gaussian design} assumption:
\begin{model}[Gaussian covariate model] 
\label{model:gcm}
In the Gaussian covariate model (GCM), we assume the inputs are independently drawn from a Gaussian distribution: 
\begin{align}
    \vec{x}^{\nu} \sim \mathcal{N}\left(\frac{\Vec{\mu}}{\sqrt{d}},\Sigma\right), \qquad \nu=1,\cdots,n.
\end{align}
We denote $\mathcal{G}_{\Vec{\theta}_0,\Sigma,\Vec{\mu}}$ the teacher-student problem under this data assumption.
\end{model}

A key motivation for this work is the common intuition that data from standard classification tasks such as MNIST are closer to multi-modal distributions than to a single-mode Gaussian. Therefore, in this manuscript we ask ourselves the question: \textit{"when is Gaussian data all you need?"}. In particular, we focus our attention to a prototypical distribution to model multi-modal data (and a universal approximator of densities): the $K$ cluster Gaussian mixture:
\begin{model}[Gaussian mixture model] 
\label{model:gmm}
In the Gaussian mixture model (GMM), we assume the inputs are independently drawn from a mixture of $K$ Gaussians:
\begin{align}
\vec{x}^{\nu}\sim \sum\limits_{c\in \mathcal{C}} p_{c}\,  \mathcal{N}\left(\frac{\vec{\mu}_{c}}{\sqrt{d}}, \Sigma_c\right), && \nu=1,\cdots, n.
    \label{eq:main:gmm}
\end{align}
where $\mathcal{C}\coloneqq \{1,\cdots, K \}$ is the set of possible clusters, and $p_{c}\in[0,1]$ is the probability of belonging to cluster $c\in\mathcal{C}$, whose means and covariance are given by $(\Vec{\mu}_c,\Sigma_c)$. 
\end{model}
We note that despite Model \ref{model:gcm} being a special case $K=1$ of Model \ref{model:gmm} when the labels $y^\nu$ do not depend on the input cluster, it will be instructive to treat the $K=1$ and $K>1$ case as two different models.

\section{Main theoretical results}
\label{sec:main:theory}
\begin{figure*}[t!]
\centering
\includegraphics[width=0.85\textwidth]{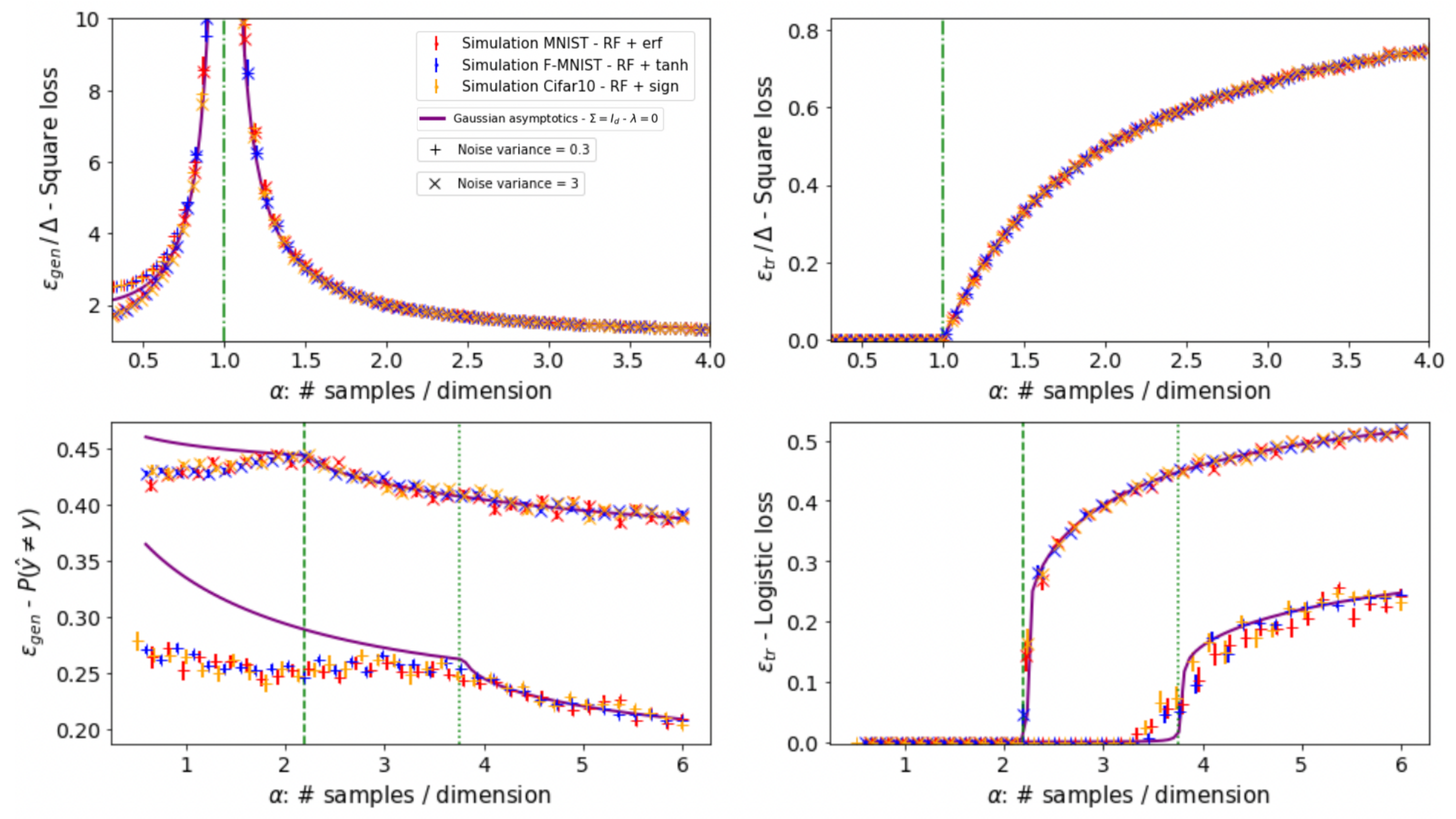}
\vspace{-1em}
    \caption{An illustration of Gaussian universality with vanishing regularization $\lambda = 0^+$ for a selection of datasets (MNIST, Fashion-MNIST, Cifar10) with a random teacher function, after a random feature map: Generalization (left) and training (right) errors as a function of the number of samples per dimension $\alpha=n/d$ for ridge regression (top pannels) and  logistic classification (bottom pannels). The solid purple line is the exact asymptotics formula for Gaussian covariates with identity covariance, while the vertical green lines are the threshold values $\alpha^\star(\Delta)$ at which an unique minimizer of the loss starts to exist. Dots show numerical simulations for different real datasets with random features maps. All data follows the Gaussian predictions for the training loss, illustrating Theorems \ref{prop:cov_universality} and \ref{prop:ridge_universality}. In particular the separability, or interpolation, threshold is the one of the Gaussian model (see corollary \ref{prop:separability}).  Each learning task was ran for two different value for the noise variance corrupting the labels - crosses are associated to $\Delta = 3$ while pluses to $\Delta = 0.3$. Error bars are built using standard deviation over 30 runs.   
  }
\label{fig:fig1}
\end{figure*}

In this section, we introduce our main theoretical results concerning universality of high-dimensional generalized linear estimation of GMMs. Our result builds on a long line of works providing an exact asymptotic characterization of empirical risk minimizers \eqref{eq:main:erm_train_err} on the proportional high-dimensional limit for Model \ref{model:gcm} \cite{dobriban_2018_high, hastie_2022_surprises, thrampoulidis_2018_precise, gcm_paper}. In particular, closer to our derivation are the rigorous results in \cite{gaussian_mix_group,gcm_paper}. We prove that the training and generalization error concentrate in high-dimensions in a deterministic expression given by the solution of a set of self-consistent equations. Then we analyze Gaussian universality, provably characterizing a set of sufficient conditions the learning task must respect such that the test and training errors of Model \ref{model:gmm} asymptotically agrees with the ones from Model \ref{model:gcm}.

Since we deal with sequences of random variables, we will need a rigorous definition of convergence. For two sequences of numbers $(a_n), (b_n)$, we write
\begin{equation}
a_n \simeq b_n \quad \text{iff} \quad \lim_{n \to \infty} a_n - b_n = 0.
\end{equation}
Accordingly, for two sequences of random variables $(X_n), (Y_n)$, we define closeness in probability by
\begin{equation}
    X_n \overset{P}{\simeq} Y_n \quad \text{iff} \quad X_n - Y_n \overset{P}{\to} 0,
\end{equation}
where $\overset{P}{\to}$ denotes convergence in probability.

\subsection{Exact asymptotics}
Our first result is to give closed-form asymptotic characterization of the performance of the minimizer of \ref{eq:main:erm_train_err} for the Gaussian Mixture model \ref{model:gmm}, generalizing the rigorous results of \cite{gaussian_mix_group}: 
\begin{proposition} (Exact asymptotics, informal statement)
\label{prop:exact_asymp_gmm}
    Consider the empirical risk minimization problem introduced in eq.~\eqref{eq:main:erm} under Gaussian mixture data given by Model \ref{model:gmm}. For any pseudo-Lispchitz performance metric $g$, the training and generalization errors \eqref{eq:main:erm_train_err} converge in the high-dimensional limit of $n,d\to\infty$ with fixed ratio $\alpha = \sfrac{n}{d}$ to deterministic expressions which are entirely determined by the solution of a set of self-consistent replica saddle-point equations \eqref{eq:app:error_formula}.
    \begin{equation}
\begin{aligned}
    \varepsilon_{\text{tr}}(\vec{\hat\theta}) \;  &\overset{P}{\simeq} \;  \varepsilon_{\text{tr}}(\vec{\theta}_0, \{\vec{\mu}_c\}_{c\in\mathcal C}, \{\Sigma_c\}_{c\in\mathcal C}) \\
    \varepsilon_{\text{gen}}(\vec{\hat\theta}) \;  &\overset{P}{\simeq} \; \varepsilon_{\text{gen}}(\vec{\theta}_0, \{\vec{\mu}_c\}_{c\in\mathcal C}, \{\Sigma_c\}_{c\in\mathcal C})
\end{aligned}
\end{equation}
\end{proposition}
The key difference between Prop. \ref{prop:exact_asymp_gmm} and Thm. 1 from \cite{gaussian_mix_group} is the distribution of the labels. While the labels in \cite{gaussian_mix_group} are given by the GMM cluster index, in Prop. \ref{prop:exact_asymp_gmm} we consider labels generated by the target function \eqref{eq:label_generation}. While we do not provide a formal proof of these results, note that the generic proof scheme from \cite{gaussian_mix_group}, that maps the solution to the study of a so-called approximate message passing algorithm \cite{donoho2009message}, can be readily adapted to our setting. Indeed, the approximate message passing scheme in our scenario is the same, and the only difference in the proof is to include a teacher in the its asymptotic analysis, similarly as was done in \cite{cornacchia2022learning}.

Once we specified the GMM, 
and properly defined the ERM in eq.~\eqref{eq:main:erm} we can evaluate the expressions in eqs.~\eqref{eq:main:erm_train_err},\eqref{eq:main:erm_gen_err} as a function of low-dimensional quantities which define the sufficient statistics of the asymptotic errors, and are also known as \textit{order parameters}. We give the detail of the computation in Appendix~\ref{sec:appendix:replica}.

\subsection{Uncorrelated teachers}
An in-depth comparison of the asymptotic expressions for the errors of Models \ref{model:gcm} \& \ref{model:gmm} reveals that their key difference lies on the way the leading target direction $\vec{\theta}_{0}$ correlate with the cluster means and covariances. A first step towards universality is therefore to characterize under which conditions the asymptotic errors are independent of the means. We make the following assumptions:
\begin{assumption}\label{assump:uncorrelated}
The teacher $\Vec{\theta}_0$ respects, $\forall (c,c') \in \mathcal{C}\times\mathcal{C}$: 
 \begin{align}
     &\lim_{n, d \to \infty} \frac{\Vec{\theta}_0^{\top}\Vec{\mu}_c}{d} = 0 \\
     &\lim_{n, d \to \infty} \frac1{d} \Vec{\theta}_0 ^{\top}\Sigma_{c'}\left(\lambda +  \sum_{c \in \mathcal C}\hat{V}^\star_{c}\Sigma_c\right)^{-1}\vec{\mu}_c \to 0
 \end{align}
 where $\{\hat{V}^{\star}_{c}\}_{c=1}^K$ are the fixed points of the (replica) saddle point equations describing the centered GMM problem.  
\end{assumption}
\begin{assumption}\label{assump:symmetric}
    The loss function, and the teacher distribution are both symmetric:
\begin{align}
    \ell(x,y) &= \ell(-x,-y) \\
    P_0\left(y|\tau\right) &=  P_0\left(-y|-\tau\right),
\end{align}
and the regularization is an $\ell^2$ penalty $\sfrac{\lambda}{2} \lVert \cdot \rVert_2^2$.
\end{assumption}

\begin{proposition}(Mean Universality)
\label{prop:mean_univers}

Under Assumptions \ref{assump:uncorrelated} and \ref{assump:symmetric},the cluster means $\{\Vec{\mu}_c\}_{c\in\mathcal{C}}$ are not relevant in high-dimensional ERM estimation:
\begin{align}
    \varepsilon^{\rm GMM}_{gen} \left(\{\Vec{\mu}_c\}_{c=1}^K,\{\Sigma_c\}_{c=1}^K\right) &\simeq \varepsilon^{\rm GMM}_{gen}\left(\Vec{0},\{\Sigma_c\}_{c=1}^K\right) \\
    \varepsilon^{\rm GMM}_{tr} \left(\{\Vec{\mu}_c\}_{c=1}^K,\{\Sigma_c\}_{c=1}^K\right) &\simeq \varepsilon^{\rm GMM}_{tr}\left(\Vec{0},\{\Sigma_c\}_{c=1}^K\right)
\end{align}
\end{proposition}
Therefore, intuitively mean universality can be achieved when the label generation process is uncorrelated with the data structure, see Appendix~\ref{sec:appendix:proofs} for a detailed discussion. The assumptions on the target and loss function are not restrictive, and are easily satisfied by odd target activation $f_0$ and margin-based losses of the form $l(y,z) = \ell(yz)$. Stronger universality can be shown by doing simplifying assumption on the data structure. Indeed, if the mixture is homogeneous (in the sense that the all the covariance are identical, a condition often called homoscedasticity in statistics) we have Gaussian universality:
\begin{theorem}{(Gaussian Universality of homoscedastic GMMs)}
\label{prop:gauss_univ}
    Under the assumption of \ref{prop:mean_univers}, consider an homoscedastic GMM:
    $$\Sigma_c = \Sigma  \qquad  \forall c \in \mathcal{C}$$
    Then for all $\alpha$ and $\Delta$, the errors of the GMM are asymptotically equal to those of a GCM:
    \begin{align}
        \varepsilon^{\rm GMM}_{\text{gen}}  \left(\{\Vec{\mu}_c\}_{c=1}^K,\Sigma \right) &\simeq 
       \varepsilon^{\rm GCM}_{\text{gen}}(\Vec{\theta}_0,\Sigma,\Vec{0}) 
        \\
        \varepsilon^{\rm GMM}_{\text{tr}} \left(\{\Vec{\mu}_c\}_{c=1}^K,\Sigma\right) &\simeq 
         \varepsilon^{\rm GCM}_{\text{tr}}(\Vec{\theta}_0,\Sigma,\Vec{0})
    \end{align}
\end{theorem}

This results follow in a straightforward way from Theorem.~\ref{prop:mean_univers}.
We have mapped under controllable assumptions a GMM problem to a simpler Gaussian one. 

Additionally, if we consider vanishing regularization, we can prove that errors are independent from the shape of the covariance, and therefore we can take an isotropic mixture:
\begin{theorem}(Covariance universality for $\lambda = 0^+$)
\label{prop:cov_universality}
Under the assumption of Thm.~\ref{prop:gauss_univ}, assume that it exists an unique minimizer of the empirical risk \eqref{eq:main:erm} with zero regularization. Then, the test \& training errors a homoscedastic GMM \eqref{model:gmm} estimation problem asymptotically coincide with those of a centered, isotropic Gaussian model $\mathcal{G}_{(\Vec{\theta}_0,\mathbf{I},\Vec{0})}$ \eqref{model:gcm}.
\end{theorem}
The proof of this general result is given in Appendix~\ref{sec:appendix:proofs}. 

Since the training error for losses such as logistic or hinge characterize the separability transition at $\lambda = 0^+$, an interesting consequence of Thm. \ref{prop:cov_universality} is the following:
\begin{corollary}(Universality of  linear separability for homoscedastic GMMs)
\label{prop:separability}
The location of the separability / interpolation transition $\alpha_{c}$ above which the data stop to be linearly separable is the same for homoscedastic GMMs and the Gaussian model.
\end{corollary}
This universality is particularly interesting in light of the detailed study of the separability transition for random teacher weights and {\it Gaussian data} in \cite{candes2020phase}. Similar universality phenomena was observed for the reconstruction transition in linear estimation in \cite{NEURIPS2019_dffbb6ef}.

Surprisingly, if we further consider a square loss minimization, the estimation of \textit{any} GMM (homoscedastic or not!) under mean universality condition can be mapped to those of a trivial  Gaussian problem: 
\begin{theorem}(Strong universality of the square loss for $\lambda =0^+$)
\label{prop:ridge_universality}
Consider underparametrized learning of GMMs with general means and covariance respecting the assumptions of Theorem.~\ref{prop:mean_univers}.  Set in eq.~\eqref{eq:main:erm}: $l(y,\hat{y}) = (y-\hat{y})^2$ and consider a generic noisy linear teacher activation $ P_0(\tau;\Delta) = \mathcal{N}(\tau,\Delta)$, with $\xi \sim \mathcal{N}(0,1)$. Then, the training error for a GMM estimation problem is given by 
\begin{align}
\label{eq:main:nolamb_ridge_closed_et}
          \varepsilon^{\rm GMM}_{tr} \left(\{\Vec{\mu}_c\}_{c=1}^K,\{\Sigma_c\}_{c=1}^K\right) &\simeq \frac{(\alpha -1)  \Delta }{\alpha}
    \end{align}
\end{theorem}
Note that the strong universality statement does not hold for the test error (a counterexample is discussed in Appendix~\ref{sec:appendix:correlated} using a strongly heteroscedastic case). However, as we see from Fig.~\ref{fig:fig1}, it seems surprisingly true that random teacher regression on real data follows the Gaussian asymptotic prediction in the underparametrized region.
Moreover, we observe in the lower panel of Fig.~\ref{fig:fig1} that even for non-quadratic losses we can draw a similar conclusion, and we investigate this further in next section.  

Finally, we note that it the limit of infinite noise and binary labels, our results give back the  ones observed for purely random Rachemacher labels proven in \cite{randlab}.

\begin{figure*}[t!]
\centering
\includegraphics[width=0.85\textwidth]{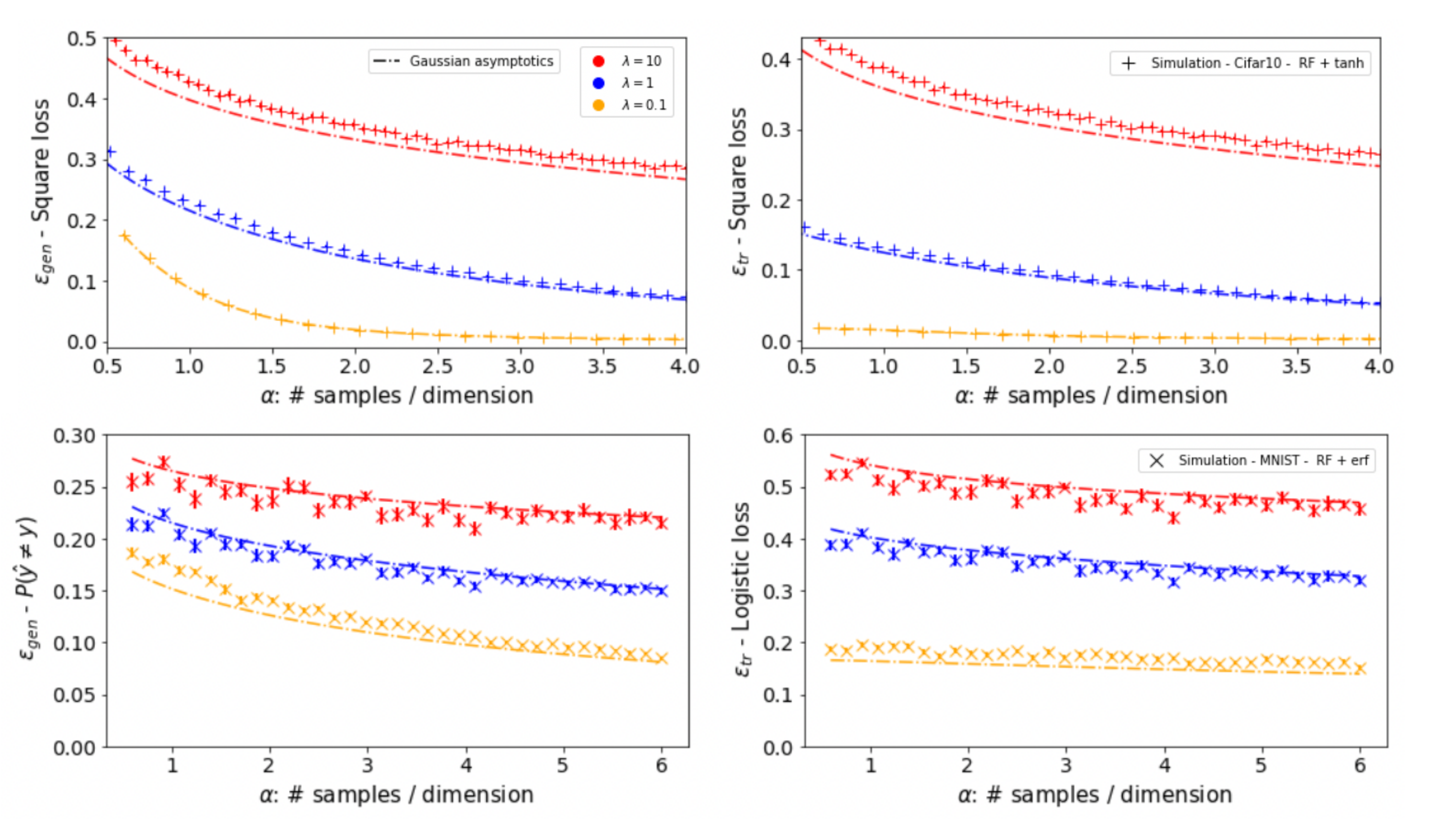}
\caption{An illustration of Gaussian universality with finite regularization $\lambda$ for a selection of datasets (MNIST, Fashion-MNIST, Cifar10) with a random teacher function, after a random feature map: Generalization (left) and training (right) errors as a function of the number of samples per dimension $\alpha=n/d$. In the upper panel we show ridge regression on Cifar10 preprocessed with RF and $\tanh$ activation, while the lower one is logistic regression on MNIST with RF and $\erf$ activation. The dashed curves are the exact asymtotics prediction of the Gaussian theory by matching covariance. The good agreement illustrates the property of theorem \ref{prop:gauss_univ}. The different colours represents different regularization strength $\lambda \in \{0.1,1,10\}$ respectively in yellow, blue, and red. Error bars are built over $30$ runs.}
\label{fig:fig2}
\end{figure*}

\subsection{Correlated teachers}

We now consider the general case where the target weights correlate with the structure in the data, relaxing the assumptions in Theorem.~\ref{prop:mean_univers}. Can we still say something? Our first result shows that the answer is yes. We focus on a simple controlled setting in which we can express the ERM performance in a closed form for \textit{any} teacher vector: 
\begin{theorem} (Exact asymptotics for isotropic covariance)
\label{prop:non_univ_mixt}
Consider a ridge regression task with a 2-clusters GMM \eqref{model:gmm}. Note that, without loss of generality we can take $\Vec{\mu}_+ = - \Vec{\mu}_- = \Vec{\mu}$. Assume isotropic covariances: 
\begin{align}
    \qquad \Sigma_+ = \Sigma_- = \mat{I}_d,
\end{align}
and denote:  
\begin{align}
\rho &= \lim\limits_{d\to\infty} \frac1d\Vec{\theta}_0^\top \Sigma \Vec{\theta}_0, \qquad \gamma = \lim\limits_{d\to\infty} \frac1d|| \Vec{\mu}||_2^2 \\ 
\pi &= \lim\limits_{d\to\infty} \frac1d\Vec{\mu}^\top \Vec{\theta}_0 \qquad P_0(\tau; \Delta) = \mathcal{N}(\tau,\Delta).
\end{align}
Defining:
\begin{align}   
      A(\eta) &= \frac{(\eta -1)^2 \left(\gamma  \eta ^2+2 \eta -1\right)}{\left(\Delta +(\eta -1)^2 \rho \right)(1+\gamma \eta)^2}\,\,\,\rm{with}\\
        \eta  &= 1 + \frac12 \left( \alpha - 1 + \lambda - \sqrt{4\lambda + \left(\alpha -1 + \lambda\right)^2} \right)\nonumber
\end{align}
The asymptotic errors admit a closed form expression in terms of $\eta$:
\begin{align}
\label{eq:pi2_scaling_eg}
    \varepsilon^{\rm GMM}_{\text{gen}} &\simeq g(\alpha,\Delta,\rho,\eta)(1 - \pi^2 A(\eta)) \\
\label{eq:pi2_scaling_et}
    \varepsilon^{\rm GMM}_{\text{tr}} &\simeq t(\alpha,\Delta,\rho,\eta)(1 - \pi^2 A(\eta)) 
\end{align}
where $g$ and $t$ are the asymptotic limit of the generalisation $\varepsilon^{\rm GCM}_{\text{gen}}(\Vec{\theta}_0,\mat{I},\Vec{0})$ and training  $\varepsilon^{\rm GCM}_{\text{tr}}(\Vec{\theta}_0,\mat{I},\Vec{0})$ error for a single Gaussian model with unit covariance and uncorrelated teacher:
\begin{align}
    g(\alpha,\Delta,\rho,\eta) & = \frac{\alpha\left(\Delta +(\eta -1)^2 \rho \right)}{\left(\alpha  - \eta ^2\right)}  \\
    t(\alpha,\Delta,\rho,\eta) &=  \frac{(\alpha-\eta)^2\left(\Delta +(\eta -1)^2 \rho \right)}{\alpha\left(\alpha - \eta^2\right) }
\end{align}
\end{theorem}
 \begin{figure*}[t!]
\includegraphics[width=1.1\textwidth,center]{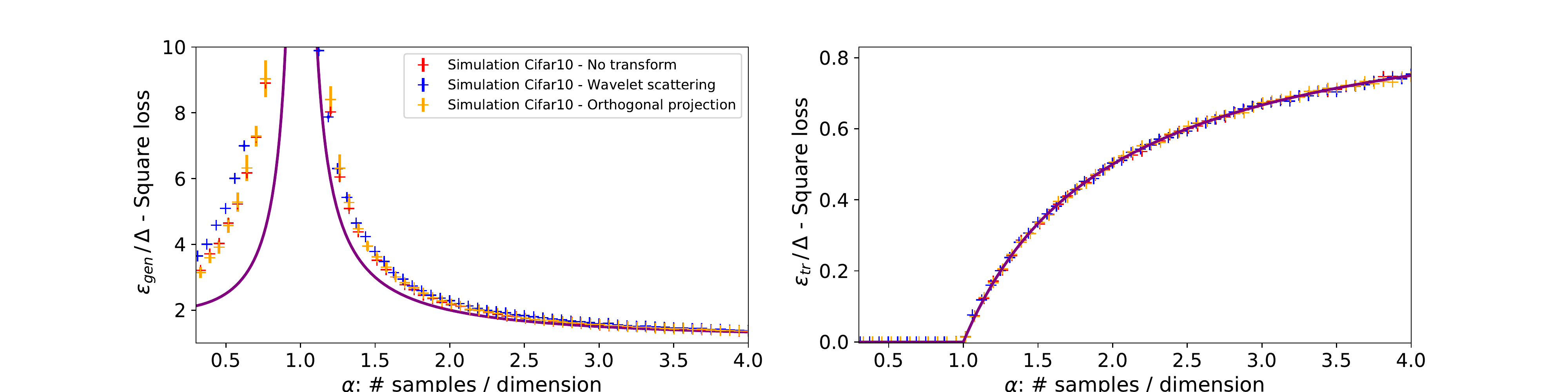}
\caption{An illustration of Gaussian universality for Training, and lack thereof for Generalization at vanishing regularization.  Instead of a random feature maps as in Fig. \ref{fig:fig1}, we pre-processed grayscale Cifar10 with wavelet scattering (blue dots) or orthogonal hadamard projections (yellow dots), or used directly the raw data (red dots). Using again a random teacher function, the training error is found to agree perfectly with the universal Gaussian asymptotics predictions with identity covariance and $\lambda = 0^+$, as for theorem \ref{prop:ridge_universality}, but the generalization is found to show clear deviation. This highlight the roles of heteroscedasticity in the data. Error bars are built over $30$ runs.}
\label{fig:rank_def}
\end{figure*}
The full closed form expressions and the extension to covariances of the form $\Sigma = \sigma \mat{I}_d$ are derived in Appendix~\ref{sec:appendix:correlated}. Note that the correction factor to the Gaussian performance scales with $\pi^2$, a measure of correlation between teacher vector and data structure. 
Hence, we would be tempted to state that targets that correlate with the data structure, i.e. $\pi \neq 0$, always break Gaussian universality. Although this is usually the case, in the limit of vanishing regularization we are in the position to present an interesting corollary of Theorem.~\ref{prop:non_univ_mixt}:
\begin{corollary}(Restoration of universality for $\lambda = 0^+$)
\label{prop:corr_teach_univers}
Consider the same setting of Theorem.~\ref{prop:non_univ_mixt}, further consider  underparametrized learning in
the limit of vanishing regularization. Then, the test and training errors for a GMM estimation are equal to a Gaussian one for any teacher vector (that is eqs (\ref{eq:pi2_scaling_eg}) and (\ref{eq:pi2_scaling_et}) at $\eta=0$):
\begin{align}
\label{eq:no_lambda_egen}
 \varepsilon^{\rm GMM}_{\text{gen}} &= \Delta \frac{\alpha}{\alpha - 1} \\
\label{eq:no_lambda_etrain}
 \varepsilon^{\rm GMM}_{\text{tr}} &= \Delta\frac{(\alpha - 1)}{\alpha}
\end{align}
\end{corollary}
We thus restore Gaussian universality for correlated teachers in the underparametrized regime, in a similar fashion to what Theorem.~\ref{prop:cov_universality} is stating for general convex losses and covariances under the mean universality condition. The universality property for correlated teachers is valid also for more general homoscedastic mixtures with identity covariance. For the sake of brevity we refer to Appendix~\ref{sec:appendix:correlated} where we discuss in detail this interesting extension. 


\section{Illustration on synthetic and real datasets }
\label{sec:main:inverstigation}
In this section we investigate the consequences of our main theoretical results. First, we consider the case of real data, illustrating the applicability of our universality theorems in cases in which the target is not correlated with the data structure. Based on these observations we move on studying simple synthetic settings described in Theorem ~\ref{prop:non_univ_mixt} and for which we can derive analytical results and systematically probe Gaussian universality. All the code used in our experiments are available in a \href{https://github.com/lucpoisson/GaussianMixtureUniversality}{GitHub repository}.
\subsection{Random teacher universality}
We analyze Gaussian universality of real data under the random teacher function assumptions of Prop.~\ref{prop:mean_univers}.  We consider three standard datasets: MNIST, Fashion-MNIST, and grayscale CIFAR-10, as well as synthetic data.
\begin{figure*}[t!]
\centering
\includegraphics[width=0.95\textwidth]{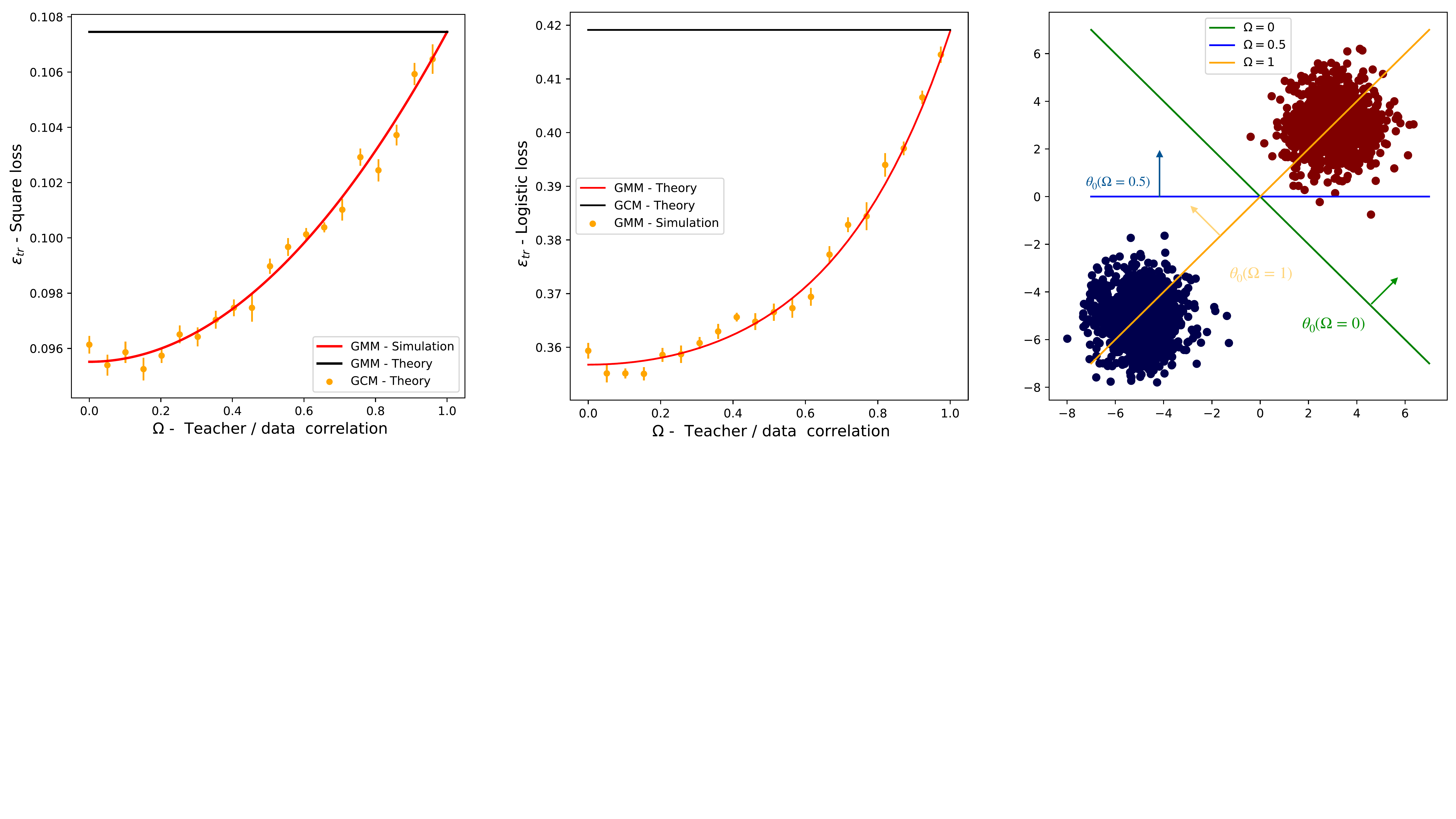}
\vspace{-12em}   
\caption{On the importance of the correlation of the task and the structure of the data: \textbf{Left \& Center }: Training errors achieved in the ERM for estimation of 2-cluster GMM with $\Vec{\mu}_+ = -\Vec{\mu}_-$ and $\Sigma_+ = \Sigma_- = \mat{I}_d$, plotted versus a correlation parameter $\Omega$: we build a series of learning tasks with teacher weights dependent on $\Omega$, namely we take $\Vec{\theta}_0(\Omega) = \Omega \Vec{\mu}_{\perp} + \sqrt{1-\Omega^2}\Vec{\mu}_{+}$, with $\Vec{\mu}_+^\top\Vec{\mu}_{\perp}=0$. We fix the dimension to be $d=500$, and $(\alpha , \lambda) = (1.2,0.7)$. The solid black line is the theoretical prediction coming from Gaussian asymptotics. The solid red line is the theoretical prediction for the GMM performance, while the orange dots are the numerical simulations which agree as expected with the theoretical prediction. In the left panel we have real labels and perform ridge regression, while on the central one we consider binary labels and perform logistic regression. The error bars are built using standard deviation over $30$ runs. \textbf{Right}: Geometrical intuition plotted for $d=2$. Three hyperplanes (lines in 2D) are displayed correspondent to $\Omega=\{0,0.5,1\}$, respectively in green, blue, and orange. As $\Omega$ increase the labels become more uncorrelated with the data structure.}
\label{fig:interpolation_exps}
\end{figure*}

\paragraph{Universality with random feature maps ---} 
To go beyond simple linear fit, we pre-process the data with a random feature map \cite{rahimi_2007_random} as follows: take an image $\Vec{\omega}_{\nu}  \in \mathbb{R}^{d^\prime}$ and map it to $\Vec{x}_{\nu} = \sigma(\mat{F}\Vec{\omega}_{\nu})$, where the elements $\mat{F} \in \mathbb{R}^{d\times d^{\prime}}$ are i.i.d drawn from $\mathcal{N}(0,1)$ and $\sigma(\cdot)$ is an activation function. As shown in \cite{rahimi_2007_random}, for $d\to\infty$ this converge to a kernel method. For each dataset we take a different non-linearity ($\erf$,$\tanh$, and $\rm{sign}$ as described in Fig.~\ref{fig:fig1}). We then build new labels by plugging in eq.~\eqref{eq:label_generation} $\Vec{\theta}_0 \sim \mathcal{N}(\Vec{0},\mat{I}_d)$, and we set for regression tasks $f_0(\tau) = \tau + \Delta$, while for classification tasks $f_0(\tau) = \rm{sign}(\tau + \Delta)$. 

First, we analyze the vanishing regularization case in Fig.~\ref{fig:fig1}: random teacher regression on pre-processed real datasets respects the strong universality of the training loss as stated in Theorem.~\ref{prop:ridge_universality}. Interestingly, we also observe a perfect match between the Gaussian asymptotics prediction and the simulations in Fig.~\ref{fig:fig1} for the generalization error. This suggests that after the random features maps the data is sufficiently close to a homoscedastic mixture. Further, the lower panel of Fig.~\ref{fig:fig1} shows that the strong universality property seems to hold beyond square loss minimization. 

We analyze as well the finite regularization setting in Fig.~\ref{fig:fig2}: we compare the simulations on real data with the prediction of the exact Gaussian asymptotics: we match the covariance of each pre-processed real dataset and compute the performance of the ERM estimator thanks to the deterministic replica formula.~\ref{prop:exact_asymp_gmm}. The predictions of the Gaussian theory  matches the numerical simulations. As discussed for Fig.~\ref{fig:fig1}, it seems that the homoscedasticity assumption in Theorem.~\ref{prop:gauss_univ} can be sometimes relaxed.

\paragraph{Universality of the double descent phenomena ---} One finding in modern machine learning that goes against the classical statistical theory wisdom is the double descent behaviour of learning curves \cite{opper1996statistical,Belkin_2019,Spigler_2019,mei19}, see upper panel of Fig.~\ref{fig:fig1}: the test error does not deteriorate as the number of parameters is increased with a characteristic divergence at $\alpha = 1$, known as the interpolation peak. As shown in Fig.~\ref{fig:fig1} and Fig.~\ref{fig:fig2}, the  generalization error and its characteristic "double descent" behavior for $\alpha>1$ is universal for homoscedastic data, while the training error appears universal even for heteroscedastic data (see Fig.~\ref{fig:rank_def}).

\paragraph{Universality of  linear separability ---}
Recently, \cite{candes2020phase} investigated the linear separability of Gaussian data with a random teacher and noise $\Delta$, generalizing the classical results by \cite{cover1965geometrical} to a single-index target. \cite{candes2020phase} has shown the existence of a critical phase transition $\alpha_c(\Delta)$ that goes continuously from $\alpha_c(\infty)=2$ (for infinite noise) to $\alpha_c>2$ for finite $\Delta$ (with $\alpha_c \to \infty$ as $\Delta \to 0$). As discussed in corollary \ref{prop:separability}, this transition is universal for homoscedastic mixtures. The interest of this transition thus extends way beyond Gaussian data and in fact agrees with the real data experiment in Figs.~\ref{fig:fig1} and \ref{fig:fig2}. As for the double descent phenomena, we believe it is a very interesting consequence of our theorems that the theoretical works mentioned above are valid way beyond the simple Gaussian assumption.

\paragraph{Non-universal behavior for generalization  ---} Finally, we implemented different transforms beyond random features for the pre-processing step: we considered the wavelet scattering transform \cite{wav_scatt}, orthogonal random projections \cite{ortho_proj}, and even no transform at all. Without the shuffling of the random projection, the fact that the data are complex and probably rather  heteroscedastic than homoscedastic, we expected a weaker form of universality. We present the results of ridge regression with labels generated by a random teacher vector in Fig.~\ref{fig:rank_def}. The strong universality statement in Theorem.~\ref{prop:ridge_universality}, that did not required any assumption on the data, remains valid and we observed a perfect collapse of the training data. However, as expected, the generalization error shows clear deviations with respect to the Gaussian behavior. Presumably these transforms do not homogenize the data enough such that Gaussian universality for the test error to hold.

\subsection{Correlated teacher} 
Previously, we showed that ridge regression with a random teacher model, even when the data is structured, can lead to universal behaviour. We now consider the dual task: take a simple homogeneous model for the data and study correlated target weights. Indeed, some recent works have studied examples of the lack of universality between Gaussian mixtures and Gaussian models \cite{tomasini_2022_failure,ingrosso22}; we want to follow this direction and use the setting described in Theorem.~\ref{prop:non_univ_mixt}: we consider a 2-cluster GMM with opposite means of norm $\Vec{\mu}$ and same covariance $\Sigma$. We build a series of learning tasks at fixed $(\alpha,\lambda)$ varying the overlap between the teacher vector and the cluster means as follows:
\begin{eqnarray}
&\Vec{\theta}_0(\Omega) = \Omega \Vec{\mu}_{\perp} + \sqrt{1-\Omega^2}\Vec{\mu} \\
&\text{such that} \qquad \Vec{\mu}^\top \Vec{\mu}_{\perp} = 0
\end{eqnarray}
The results are presented in Fig.~\ref{fig:interpolation_exps} for both ridge and logistic regression. They clearly show that a small correlation between the target weights and the mixture means breaks Gaussian Universality. A neat geometrical intuition of the result is given in Fig.~\ref{fig:interpolation_exps}: as we decrease $\Omega$, we generate labels which are more and more correlated with the data distribution, and consequently there is a correction factor to the Gaussian prediction as  Theorem.~\ref{prop:non_univ_mixt} predicts.    
We refer to Appendix.~\ref{sec:appendix:correlated} for a more detailed discussion on correlated teachers. We conclude that universality is broken in tasks where the labels are correlated with the structure. Note, however, that as proven in Corollary.~\ref{prop:corr_teach_univers}: the discrepancy between the GMM prediction and the Gaussian one goes to zero \textit{for any} correlation measure $\Omega$ as $\lambda \to 0^+$, where the universality is restored as suggested by Thm.~\ref{prop:ridge_universality}.

\section{Acknowledgments}
\label{sec:main:acknowledgement}
We thank Yatin Dandi, Federica Gerace, Sebastian Goldt, Denny Wu \& Lenka Zdeborov\'a for useful discussions. We acknowledge funding from the Swiss National Science Foundation grant SNFS OperaGOST, $200021\_200390$ and the \textit{Choose France - CNRS AI Rising Talents} program.
\pagebreak

\newpage 
\bibliography{example_paper}
\bibliographystyle{unsrt}

\newpage
\appendix
\onecolumn
\section{Replica computation}
\label{sec:appendix:replica}
\subsection{Formal statement of the theorem}\label{subsec:appp:formal_theorem}
We first provide the full statement of Proposition \ref{prop:exact_asymp_gmm}. Consider a minimization problem of the form
\begin{equation} 
\vec{\hat\theta} = \argmin_{\vec{\theta \in \mathbb{R}^d}} \frac{1}{n}\sum\limits_{\nu=1}^{n}\ell\left(y^{\nu}, \vec{\theta}^{\top}\vec{x}^\nu\right) + r(\vec{\theta}),
\end{equation}
where the data $(\vec{x}^\nu, y^\nu)$ is generated according to the following Gaussian mixture model:
\begin{equation}
    \vec{x}^\nu \sim \sum_{c\in\mathcal C} p_c\, \mathcal N\left( \frac{\vec{\mu}_c}{\sqrt{d}}, \Sigma_c \right) \quad \text{and} \quad y^\nu \sim P_0(\cdot \,|\, \vec{\theta}_0^\top \vec{x}^\nu)
\end{equation}
and assume the following:
\begin{enumerate}
    \item The functions $\ell(y, \cdot)$ and $r$ are continuous and coercive, and the function $\ell(y, \cdot) + r(\cdot)$ is strongly convex,
    \item The covariance matrices $\Sigma_c$ are positive definite, and their spectral norms are uniformly bounded,
    \item The means $\vec{\mu_c} / \sqrt{d}$ and the teacher vector $\vec{\theta}_0$ are uniformly bounded,
    \item the number of clusters $|\mathcal C|$ is finite,
    \item the distribution $P_0$ is sub-gaussian with uniformly bounded norm.
\end{enumerate}
Then, as $n, d \to \infty$ with $n/d \to \alpha > 0$, we have
\begin{equation}\label{eq:app:error_formula}
\begin{aligned}
    \varepsilon_{\text{tr}}(\vec{\hat\theta}) \;  &\overset{P}{\simeq} \; \sum_{c \in \mathcal C} p_c \, \mathbb E_{\omega^{(s)}_c, \omega^{(t)}_c, y}\left[ \ell\left(y, \prox_{V_c^\star\ell(y, \cdot)}(\omega^{(s)}_c) \right)\right] &=: \varepsilon_{\text{tr}}(\vec{\theta}_0, \{\vec{\mu}_c\}_{c\in\mathcal C}, \{\Sigma_c\}_{c\in\mathcal C}) \\
    \varepsilon_{\text{tr}}(\vec{\hat\theta}) \;  &\overset{P}{\simeq} \; \sum_{c \in \mathcal C} p_c \, \mathbb E_{\omega^{(s)}_c, \omega^{(t)}_c, y}\left[ \ell\left(y, \omega^{(s)}_c\right)\right] &=: \varepsilon_{\text{tr}}(\vec{\theta}_0, \{\vec{\mu}_c\}_{c\in\mathcal C}, \{\Sigma_c\}_{c\in\mathcal C})
\end{aligned}
\end{equation}
where
\begin{equation}
    \begin{pmatrix}
    \omega_c^{(s)} \\ \omega_c^{(t)}
    \end{pmatrix} \sim \mathcal N \left( \begin{bmatrix}
        \pi_c \\ h^\star_c
    \end{bmatrix}, \begin{bmatrix}
        \rho & m^\star_c \\
        m^\star_c & q^\star_c
    \end{bmatrix}\right) \quad \text{and} \quad y \sim P_0\left( \cdot \,\big|\, \omega_c^{(t)}\right).
\end{equation}
and the $\prox$ function is
\begin{equation}
    \prox_{f}(x) = \min_{z \in \mathbb R} \left[\frac12 (z-x)^2 + f(z) \right]
\end{equation}

The overlaps used in the equation are defined as follows:
\begin{equation}
    \rho_c = \frac1d \vec{\theta}_0^\top \Sigma_c \vec{\theta}_0, \quad \pi_c = \frac1d \vec{\theta}_0^\top \vec{\mu}_c,
\end{equation}
and $(V^{\star}_{c},q^{\star}_{c},m^{\star}_{c},h^{\star}_{c})_{c\in\mathcal{C}}$ are the unique fixed point of the following set of self-consistent \emph{replica saddle-point equations}:
\begin{equation}\label{eq:saddle_point}
\begin{split}
&\begin{cases}
	\hat{V}_{c} &= \alpha p_{c} \mathbb{E}_{\xi_c\sim\mathcal{N}(0,1)}\left[\int\dd y~\mathcal{Z}_{0}\left(y,\pi_{c} + \frac{m_{c}}{\sqrt{q}_{c}} \xi_c,\rho-\frac{m_{c}^2}{q_{c}}\right)\partial_{\omega}f_{\ell}(y,h_{c}+\sqrt{q_{c}}\xi_c,V_{c})\right]\\
	\hat{q}_{c} &= \alpha p_{c} \mathbb{E}_{\xi_c\sim\mathcal{N}(0,1)}\left[\int\dd y~\mathcal{Z}_{0}\left(y,\pi_{c} + \frac{m_{c}}{\sqrt{q}_{c}} \xi_c,\rho-\frac{m_{c}^2}{q_{c}}\right) f_{\ell}(y,h_{c}+\sqrt{q_{c}}\xi_c,V_{c})^2\right]\\
	\hat{m}_{c} &= \alpha p_{c} \mathbb{E}_{\xi_c\sim\mathcal{N}(0,1)}\left[\int\dd y~ 	\partial_{\omega}\mathcal{Z}_{0}\left(y,\pi_{c} + \frac{m_{c}}{\sqrt{q}_{c}} \xi_c,\rho-\frac{m_{c}^2}{q_{c}}\right) f_{\ell}(y,h_{c}+\sqrt{q_{c}}\xi_c,V_{c})\right]  \\
	\hat{h}_{c} &= \alpha p_{c} \mathbb{E}_{\xi_c\sim\mathcal{N}(0,1)}\left[\int\dd y~ 	\mathcal{Z}_{0}\left(y,\pi_{c} + \frac{m_{c}}{\sqrt{q}_{c}} \xi_c,\rho-\frac{m_{c}^2}{q_{c}}\right) f_{\ell}(y,h_{c}+\sqrt{q_{c}}\xi_c,V_{c})\right] 
\end{cases}	\\[1em]
&\begin{cases}
	V_{c} &= \mathbb{E}_{\{\vec{\xi}_c\}~\stackrel{\text{i.i.d}}{\sim}~\mathcal{N}(\vec{0},\mat{I}_{d})}\left[\vec{\hat\eta}^\top \hat q_c^{-1/2}\Sigma_c^{1/2} \vec{\xi}_c\right]\\
	q_{c} &= \mathbb{E}_{\{\vec{\xi}_c\}~\stackrel{\text{i.i.d}}{\sim}~\mathcal{N}(\vec{0},\mat{I}_{d})}\left[\vec{\hat\eta}^\top \Sigma_c \vec{\hat\eta}\right]\\
	m_{c} &= \mathbb{E}_{\{\vec{\xi}_c\}~\stackrel{\text{i.i.d}}{\sim}~\mathcal{N}(\vec{0},\mat{I}_{d})}\left[\vec{\theta}_{0}^{\top}\Sigma_{c} \vec{\hat\eta}\right]   \\
	h_c &= \mathbb{E}_{\{\vec{\xi}_c\}~\stackrel{\text{i.i.d}}{\sim}~\mathcal{N}(\vec{0},\mat{I}_{d})}\left[\vec{\mu}_{c}^{\top}\vec{\hat\eta}\right]
\end{cases}	
\end{split}
\end{equation}
and we have defined the following auxiliary functions:
\begin{align}
        \mathcal{Z}_{0}(y, \omega, V) &= \int \frac{\dd\lambda}{\sqrt{2\pi V}} P_{0}(y|\lambda)~e^{-\frac{(\lambda-\omega)^2}{2V}}\\
        f_{\ell}(y,\omega, V)  &= -V^{-1}  \partial_\omega \mathcal{M}_{V\ell(y, \cdot)}(\omega) = V^{-1}\left( \prox_{V\ell(y, \cdot)}(\omega) - \omega \right).
\end{align}
Finally, the auxiliary variable $\vec{\hat \eta}$ is a function of $\{\vec{\xi}_c, \hat V_c, \hat m_c, \hat q_c, \hat h_c\}$:
\begin{equation}\label{eq:app:def_etahat}
    \vec{\hat \eta} = \prox_{r(\hat \Sigma^{-1/2} \cdot)} \left( \hat \Sigma^{-1/2} \left( \sum_{c \in \mathcal C} \vec{\hat \mu}_c + \sqrt{\hat q_c} \Sigma_c^{1/2} \vec{\xi}_c \right) \right)
\end{equation}
with
\begin{equation}\label{eq:app:def_muhat_sigmahat}
    \vec{\hat \mu_c} = \hat h_c \vec{\mu_c} + \hat m_c \Sigma_c \vec{\theta}_0 \quad \text{and} \quad \hat \Sigma = \sum_{c \in \mathcal C} \hat V_c \Sigma_c
\end{equation}

We will not prove this theorem, since the proof is virtually equivalent to the one in \cite{gaussian_mix_group, cornacchia2022learning}. Instead, the next sections are dedicated to the derivation of these equations, using the so-called \emph{replica method} from statistical physics \cite{mezard.parisi.ea_1987_spin}. Simpler particular cases of these equations can be found in Appendix \ref{subsec:app:examples}.

\subsection{Gibbs measure and free energy}
The starting point for our replica computation is to define the Gibbs measure over the weights $\mat{W}\in\mathbb{R}^{K\times d}$:
\begin{align}
    \mu_{\beta}(\vec{\theta}) = \frac{1}{\mathcal{Z}_{\beta}}e^{-\beta \left[\sum\limits_{\nu=1}^{n}\ell\left(y^{\nu}, \vec{\theta}^{\top}\vec{x}^{\nu}\right) + \lambda r(\vec{\theta})\right]} = \frac{1}{\mathcal{Z}_{\beta}}e^{-\beta r(\vec{\theta})}\prod\limits_{\nu=1}^{n} e^{-\beta \ell(y^{\nu}, \vec{\theta}^{\top}\vec{x}^{\nu})}
\end{align}
\noindent where $\beta>0$ is a parameter we eventually want to send to $\beta\to 0^{+}$, and $\mathcal{Z}_{\beta}\in\mathbb{R}$ is the partition function (the normalisation of the Gibbs measure). For convenience, we will define the following useful notation:
\begin{align}
    P_{\vec{\theta}}(\vec{\theta}) \equiv e^{-\beta r(\vec{\theta})}, && P_{\ell}(y^{\nu}|\vec{\theta}^{\top}\vec{x}^{\nu}) \equiv e^{-\beta \ell(y^{\nu}, \vec{\theta}^{\top}\vec{x}^{\nu})}
\end{align}
which allow us to write the Gibbs measure as:
\begin{align}
    \mu_{\beta}(\vec{\theta}) = \frac{1}{\mathcal{Z}_{\beta}}P_{\vec{\theta}}(\vec{\theta})\prod\limits_{\nu=1}^{n} P_{\ell}(y^{\nu}|\vec{\theta}^{\top}\vec{x}^{\nu})
\end{align}
When $\beta \to \infty$, the Gibbs measure $\mu_\beta$ will concentrate around the $\vec{\theta}$ that minimize the empirical risk \eqref{eq:main:erm}. As a result, the free energy density defined as
\begin{align}\label{eq:app:def_free_energy}
  -\beta f_{\beta} = \lim\limits_{d\to\infty}\frac{1}{d}\mathbb{E}\log\mathcal{Z}_{\beta}
\end{align}
\noindent where the limit is taken with $n/d = \alpha$ fixed, and the expectation is over the distribution of the data, will concentrate around the minimum risk. In order to take the expectation explicitly, we use the replica trick:
\begin{align}
    \log{\mathcal{Z}_{\beta}} = \lim\limits_{s\to 0^{+}}\frac{1}{s}\partial_{s}\mathcal{Z}_{\beta}^{s}|_{s=0}
\end{align}
Therefore, the replica computation boils down to the computation of the averaged replicated partition function:
\begin{align}
    \mathbb{E}\mathcal{Z}_{\beta}^{s} &= \mathbb{E}\left[\int_{\mathbb{R}^{d}}\dd\vec{\theta}~P_{\vec{\theta}}(\vec{\theta}) \prod_{\nu = 1 }^nP_{\ell}\left(y^{\nu}|\vec{\theta}^{\top}\vec{x}^{\nu}\right)\right]^{s}\\
    &=\prod\limits_{\nu=1}^{n}\mathbb{E}_{(\vec{x}^{\nu},y^{\nu})}\int\prod\limits_{a=1}^{s}\dd\vec{\theta}^{a}~P_{\vec{\theta}}(\vec{\theta}^{a}) P_{\ell}\left(y^{\nu}|{\vec{\theta}^{a}}^{\top}\vec{x}^{\nu}\right)
    \label{eq:app:def_replica_partition}
\end{align}

\subsection{Taking the average over the data}
The main difference in this computation with respect to the usual binary Gaussian mixture is that the labels are generated by a teacher target function. In other words, the joint distribution of the data is given by:
\begin{align}
	P_{\vec{\theta}_{0}}(\vec{x},y) = P_0(y|\tau)\sum\limits_{c\in \mathcal{C}} \rho_{c} \, \mathcal{N}\left(\vec{x}; \frac{\vec{\mu}_{c}}{\sqrt{d}}, \mat{I}_{d}\right)
\end{align}
\noindent where $P_0(y|\tau)$ is the probability induced by the teacher activation $f_0(\tau)$.
The expression in \eqref{eq:app:def_replica_partition} can then be written as:
\begin{align}
    \mathbb{E}\mathcal{Z}_{\beta}^{s} &=\int\prod\limits_{a=1}^{s}\dd\vec{\theta}^{a}~P_{\vec{\theta}}(\vec{\theta}^{a}) \left(\int_{\mathbb{R}} \dd{y}\int_{\mathbb{R}^d} \dd{\vec{x}} P_{0}(\vec{y}|\vec{\theta}_{0}^{\top}\vec{x})P_{\vec{x}}(\vec{x})  \prod\limits_{a=1}^{s}P_{\ell}\left(y^{\nu}|{\vec{\theta}^{a}}^{\top}\vec{x}^{\nu}
    \right) \right)^n
\end{align}
We make a change of variables by introducing the \emph{local fields} $\tau = \vec{\theta}_0^\top \vec{x}$ and $\lambda^a = \vec{\theta}^{a\top} \vec{x}$. The distribution of the local fields will be a Gaussian mixture itself, and we can compute the moments for each mode:
\begin{equation}
\begin{split}
    \pi_{c} &\equiv \mathbb{E}_{\vec{x}^{\nu}}\left[\tau_c\right] = \frac{1}{d}\vec{\theta}_{0}^{\top}\vec{\mu}_{c}\\ 
    \rho_{c} &\equiv \text{Var}_{\vec{x}^{\nu}}\left[\tau_c\right] = \frac{1}{d}\vec{\theta}_{0}^{\top}\Sigma_{c}\vec{\theta}_{0}\\
    h^{a}_{c} &\equiv \mathbb{E}_{\vec{x}^{\nu}}\left[\lambda^{a}_c\right] = \frac{1}{d}{\vec{\theta}^{a}}^{\top}\vec{\mu}_{c}\\ 
    q_{c}^{ab} &\equiv \text{Cov}_{\vec{x}^{\nu}}\left[\lambda_c^{a},\lambda_c^{b}\right] = \frac{1}{d}{\vec{\theta}^{a}}^{\top}\Sigma_{c}\vec{\theta}^{b}\\
    m_{c}^{a} &\equiv \text{Cov}_{\vec{x}^{\nu}}\left[\tau_c,\lambda_c^{a}\right] = \frac{1}{d}\vec{\theta}_{0}^{\top}\Sigma_{c}\vec{\theta}^{a}
    \end{split}
\end{equation}
Equivalently, the local fields distribution is the following low-dimensional Gaussian mixture distribution:
\begin{align}
    P(\tau, \vec{\lambda}) = \sum\limits_{c\in\mathcal{C}}p_{c}~\mathcal{N}\left(\begin{bmatrix}\pi_{c}\\\vec{h}_{c}\end{bmatrix}, \begin{bmatrix}\rho_{c} & \vec{m}_{c}^{\top}\\ \vec{m}_{c} & \mat{Q}_{c}\end{bmatrix}\right),
    \label{eq:local_field_distr}
\end{align}
\noindent where we have defined the vectors $\vec{m}, \vec{h} \in\mathbb{R}^{s}$ with entries $m^{a}$, $h^{a}$ and the matrix $\mat{Q}\in\mathbb{R}^{s\times s}$ with entries $q^{ab}$.
Therefore, we can factorize the partition function by only integrating over the local field distribution:
\begin{align}
    \mathbb{E}\mathcal{Z}_{\beta}^{s} &=\int\prod\limits_{a=1}^{s}\dd\vec{\theta}^{a}~P_{\vec{\theta}}(\vec{\theta}^{a})\left[\int\dd y\int\dd\tau~P_{0}(y|\tau)\int\left(\prod\limits_{a=1}^{s}\dd\lambda^{a} P_{\ell}(y|\lambda^{a}) \right)P(\tau,\vec{\lambda}|\vec{m},\vec{h},\mat{q})\right]^{n}
\end{align}

\subsection{Writing as a saddle-point problem}
Note that $(\pi_{c},\rho_{c})$ are fixed inputs in the problem. By Fourier transform arguments, we can write
\begin{equation}
    \delta\left( m_c^a - \frac{1}{d}\vec{\theta}_{0}^{\top}\Sigma_{c}\vec{\theta}^{a} \right) = \int_{i\mathbb{R}} \frac{\dd{\hat m_{c}^a}}{2\pi} \, e^{\hat m_{c}^a \left( dm_c^a - \vec{\theta}_{0}^{\top}\Sigma_{c}\vec{\theta}^{a} \right)},
\end{equation}
where the integral is on the imaginary line. By doing the same arguments for the $q_c^{ab}$ and $h_c^a$, it ensues that
\begin{align}
     \mathbb{E}\mathcal{Z}_{\beta}^{s} &= \int\prod_{c \in \mathcal C}\prod\limits_{a=1}^{s}\frac{\dd m_{c}^{a}\hat{m}_{c}^{a}}{2\pi}\frac{\dd h_{c}^{a}\hat{h}_{c}^{a}}{2\pi}\prod\limits_{1\leq a\leq b\leq s}\frac{\dd q_{c}^{ab}\dd\hat{q}_{c}^{ab}}{2\pi}e^{d\Phi^{(s)}_{\rm{rs}}(\vec{m}, \hat{\vec{m}}, \mat{q}, \hat{\mat{q}})}
\end{align}
\noindent where we have defined the free energy potential:
\begin{align}
    \Phi^{(s)}_{\rm{rs}}(\vec{m}, \hat{\vec{m}}, \mat{q}, \hat{\mat{q}}) &=  \sum\limits_{c\in\mathcal{C}}\left[\sum\limits_{a=1}^{s}\hat{m}_{c}^{a}m_{c}^{a} + \sum\limits_{a=1}^{s}\hat{h}_{c}^{a}h_{c}^{a}+ \sum\limits_{1\leq a\leq b\leq s}\hat{q}_{c}^{ab}q_{c}^{ab}\right] - \alpha \Psi^{(s)}_{g}(\vec{m}, \vec{h}, \mat{q}) - \Psi^{(s)}_{\theta}(\hat{m},\hat{h},\hat{q})\notag\\
    \Psi^{(s)}_{\ell}(\vec{m},\vec{h},\mat{q}) &\equiv \log\int\dd y\int\dd\tau~P_{0}(y|\tau)\int\prod\limits_{a=1}^{s}\dd\lambda^{a} P_{\ell}(y|\lambda^{a}) P(\tau,\vec{\lambda}|\vec{m},\vec{h},\mat{q})\notag\\
    \Psi^{(s)}_{\theta}(\hat{\vec{m}},\hat{\vec{h}}, \hat{\mat{q}}) &\equiv \frac{1}{d}\log \int\prod\limits_{a=1}^{s}\dd\vec{\theta}^{a}~P_{\vec{\theta}}(\vec{\theta})\prod\limits_{c\in\mathcal{C}} e^{\sum\limits_{a=1}^{s}\left[\hat{h}_{c}^{a} {\vec{\theta}^{a}}^{\top}\vec{\mu}_{c}+\hat{m}^{a}\vec{\theta}_{0}^{\top}\Sigma_{c}\vec{\theta}^{a}\right] + \sum\limits_{1\leq a\leq b\leq s}\hat{q}_{c}^{ab}{\vec{\theta}^{a}}^{\top}\Sigma_{c}\vec{\theta}^{b}}
\end{align}
Therefore, as we take the $d\to\infty$ limit, we can apply the saddle-point method \cite{mezard.parisi.ea_1987_spin} to compute $f_\beta$
\begin{align}
    -\beta f_{\beta} = \underset{\vec{m}, \hat{\vec{m}}, \vec{h},\hat{\vec{h}}, \mat{q},\hat{\mat{q}}}{\extr} ~\lim\limits_{s\to 0^{+}}\frac{\Phi^{(s)}(\vec{m}, \hat{\vec{m}}, \vec{h},\hat{\vec{h}}, \mat{q}, \hat{\mat{q}})}s
\end{align}

\begin{remark}
We have introduced the parameters $\hat q, \hat m, \hat h$ as pure imaginary numbers, but the optimization problem considers them as real numbers. This stems from the fact that the function $\Phi$ is holomorphic, so the integral is independent from the contour of integration. More details can be found in \cite{talagrand_2022_what}.
\end{remark}

\subsection{Replica symmetric ansatz}
In order to make progress with the $s\to 0^{+}$ limit, we make the following replica symmetric ansatz:
\begin{align}
    m_{c}^{a} = m_{c} && \hat{m}^{a}_{c} = \hat{m}_{c}, && a =1, \cdots, s\notag\\
    h_{c}^{a} = h_{c} && \hat{h}^{a}_{c} = \hat{h}_{c}, && a =1, \cdots, s\notag\\
    q_{c}^{aa} = r_{c}, && \hat{q}^{aa}_{c} = -\frac{1}{2}\hat{r}_{c}, && a=1,\cdots, s \notag\\
    q_{c}^{ab} = q_{c}, && \hat{q}^{aa}_{c} = \hat{q}_{c}, && 1\leq a < b\leq s 
\end{align}
for all $c\in\mathcal{C}$. Since $\ell$ and $r$ are convex, the function $\Phi^{(s)}$ has a unique saddle point, which must therefore coincide with the replica symmetric ansatz. We also define
\begin{align}
    V_c = r_c - q_c, \qquad \text{and}  \qquad \hat{V}_c = \hat{r}_c - \hat{q}_c
\end{align}
By inserting this ansatz above, we can take the $s \to 0^+$ limit. This is done through a classical but computationally heavy method known as the \emph{Hubbard-Stratonovich transform}; details can be found in \cite{gerace20a}, Appendix C. We simply reproduce the final results here:
\begin{align}
-\beta f_{\beta} = \underset{\{m_c, \hat{m_c}, h_c,\hat{h}_c, \mat{q}_c,\hat{\mat{q}_c}, V_c , \hat{V}_c\}}{\text{extr}} ~\Phi_{\text{rs}}(\{m_c, \hat{m_c}, h_c,\hat{h}_c, \mat{q}_c,\hat{\mat{q}_c}, V_c , \hat{V}_c\})
\end{align}
The replicated free energy $\Phi_{\text{rs}}$ is given by the following formula:
\begin{align}
    \Phi_{\text{rs}}(\{m_c, \hat{m_c}, h_c,\hat{h}_c, \mat{q}_c,\hat{\mat{q}_c}, V_c , \hat{V}_c\}) = &\sum\limits_{c\in\mathcal{C}}\left[\frac{1}{2}\left( \hat{V}_{c}q_{c} - \hat{q}_{c}V_{c}\right) -\frac{1}{2}\hat{V}_{c}V_{c} + \hat{m}_{c}m_{c} +\hat{h}_c h_c\right] \\
    &- \alpha \Psi_{\ell}(\{m_c, h_c, q_c, V_c\}) - \Psi_{\vec{\theta}}(\{\hat{m_c},\hat{h}_c,\hat{q_c}, \hat{V}_c\})\notag\\
\end{align}
where we have decomposed the contributions coming from the loss ($\Psi_\ell$) and the regularization ($\Psi_{\vec{\theta}}$):
\begin{align}
\label{eq:app:def_free_loss}
   \Psi_{\ell} &= \sum\limits_{c\in\mathcal{C}}\mathbb{E}_{\xi_c~\stackrel{\text{i.i.d}}{\sim}~\mathcal{N}(0,1)}\int\dd y\mathcal{Z}_{0}\left(y, \frac{m_c}{\sqrt{q_c}}\xi+\pi_{c}, \rho-\frac{m_{c}^{2}}{q_{c}}\right)\log\mathcal{Z}_\ell(y,h_{c}+\sqrt{q_{c}}\,\xi_c,V_{c})\\
    \Psi_{\vec{\theta}} &= \frac{1}{d}\,\mathbb{E}_{\vec{\xi}_c~\stackrel{\text{i.i.d}}{\sim}~\mathcal{N}(\vec{0},\mat{I}_{d})}\log \mathcal{Z}_{\vec{\theta}}\left(\sum\limits_{c\in\mathcal{C}}\hat{V}_{c}\Sigma_{c}, \sum\limits_{c\in\mathcal{C}}\hat{h}_{c}\vec{\mu}_{c}+\hat{m}_{c}\Sigma_{c}\vec{\theta}_0+\sqrt{\hat{q}_{c}}\Sigma_{c}^{1/2}\vec{\xi}_c\right),
\end{align}
and defined the following auxiliary free energies:
\begin{align}
    \mathcal{Z}_{\ell/0}(y, \omega, V) &= \int \frac{\dd\lambda}{\sqrt{2\pi V}} P_{\ell/0}(y|\lambda)~e^{-\frac{(\lambda-\omega)^2}{2V}}\\
     \mathcal{Z}_{\vec{\theta}}(\mat{A},\vec{b}) &= \int \dd\vec{\theta} e^{-\beta r(\vec{\theta})} e^{-\frac{1}{2}\vec{\theta}^{\top}\mat{A}\vec{\theta}+\vec{b}^{\top}\vec{\theta}}
\end{align}

\subsection{Taking the zero temperature limit}
In order to take the $\beta \to \infty$ limit, we make the following rescalings:
\begin{align}\label{eq:app:beta_rescaling}
V_c&\to \beta^{-1} V_c & \hat{V}_c&\to \beta\hat{V}_c & \hat{q}_c&\to \beta^2\hat{q}_c & \hat{m}_c&\to \beta\hat{m}_c.
\end{align}
It is easy to check how this rescaling affects $\beta^{-1}\Phi^{(\text{rs})}$, so we only need to consider $\Psi_\ell$ and $\Psi_{\vec{\theta}}$. 

We start with the latter: letting
\[ \mathcal L_{\vec{\theta}}(\vec{\theta}) = \frac12\, \vec{\theta}^\top A \vec{\theta} - \vec{b}^\top \vec{\theta} + r(\vec{\theta})  \]
by the Laplace method for any $A, \vec{b}$,
\begin{align}
\lim\limits_{\beta\to\infty}\frac{1}{\beta}	\log Z_{\vec{\theta}}(\beta A, \beta \vec{b}) = -\inf_{\vec{\theta}} \mathcal L_{\vec{\theta}}(\vec{\theta}) = -\mathcal L_{\vec{\theta}}(\vec{\eta})
\end{align}
\noindent where
\begin{equation}
    \vec{\eta} = \prox_{r(A^{1/2}\cdot)}(A^{1/2} \vec{b}).
\end{equation}
As a result, every integral involved in $\Psi_{\vec{\theta}}$ (and, later, its partial derivatives) concentrates around its value at $\vec{\hat \eta}$ defined in \eqref{eq:app:def_etahat}.
For the term $\Psi_\ell$, the term $\mathcal Z_0$ is left unchanged, and by the same reasoning as above
\begin{equation}
\lim\limits_{\beta\to\infty}\frac{1}{\beta}	\log\mathcal{Z}_{\ell}(y,\omega,\beta V) = -V^{-1}\mathcal{M}_{V\ell(y,\cdot)}(\omega).
\end{equation}

\subsection{Saddle-point equations}
The saddle-point equations are obtained by taking the derivatives of the free energy potential with respect to the overlap parameters. We obtain a set of self-consistent equations which we should solve in order to find a fixed point: 
\begin{align}
\label{eq:saddle:hat}
&\begin{cases}
	\hat{V}_{c} &= \alpha p_{c} \mathbb{E}_{\xi_c\sim\mathcal{N}(0,1)}\left[\int\dd y~\mathcal{Z}_{0}\left(y,\pi_{c} + \frac{m_{c}}{\sqrt{q}_{c}} \xi_c,\rho-\frac{m_{c}^2}{q_{c}}\right)\partial_{\omega}f_{\ell}(y,h_{c}+\sqrt{q_{c}}\xi_c,V_{c})\right]\\
	\hat{q}_{c} &= \alpha p_{c} \mathbb{E}_{\xi_c\sim\mathcal{N}(0,1)}\left[\int\dd y~\mathcal{Z}_{0}\left(y,\pi_{c} + \frac{m_{c}}{\sqrt{q}_{c}} \xi_c,\rho-\frac{m_{c}^2}{q_{c}}\right) f_{\ell}(y,h_{c}+\sqrt{q_{c}}\xi_c,V_{c})^2\right]\\
	\hat{m}_{c} &= \alpha p_{c} \mathbb{E}_{\xi_c\sim\mathcal{N}(0,1)}\left[\int\dd y~ 	\partial_{\omega}\mathcal{Z}_{0}\left(y,\pi_{c} + \frac{m_{c}}{\sqrt{q}_{c}} \xi_c,\rho-\frac{m_{c}^2}{q_{c}}\right) f_{\ell}(y,h_{c}+\sqrt{q_{c}}\xi_c,V_{c})\right]  \\
	\hat{h}_{c} &= \alpha p_{c} \mathbb{E}_{\xi_c\sim\mathcal{N}(0,1)}\left[\int\dd y~ 	\mathcal{Z}_{0}\left(y,\pi_{c} + \frac{m_{c}}{\sqrt{q}_{c}} \xi_c,\rho-\frac{m_{c}^2}{q_{c}}\right) f_{\ell}(y,h_{c}+\sqrt{q_{c}}\xi_c,V_{c})\right] 
\end{cases}	\\[1em]
&\begin{cases}
	V_{c} &= \mathbb{E}_{\{\vec{\xi}_c\}~\stackrel{\text{i.i.d}}{\sim}~\mathcal{N}(\vec{0},\mat{I}_{d})}\left[\vec{\hat\eta}^\top \hat q_c^{-1/2}\Sigma_c^{1/2} \vec{\xi}_c\right]\\
	q_{c} &= \mathbb{E}_{\{\vec{\xi}_c\}~\stackrel{\text{i.i.d}}{\sim}~\mathcal{N}(\vec{0},\mat{I}_{d})}\left[\vec{\hat\eta}^\top \Sigma_c \vec{\hat\eta}\right]\\
	m_{c} &= \mathbb{E}_{\{\vec{\xi}_c\}~\stackrel{\text{i.i.d}}{\sim}~\mathcal{N}(\vec{0},\mat{I}_{d})}\left[\vec{\theta}_{0}^{\top}\Sigma_{c} \vec{\hat\eta}\right]   \\
	h_c &= \mathbb{E}_{\{\vec{\xi}_c\}~\stackrel{\text{i.i.d}}{\sim}~\mathcal{N}(\vec{0},\mat{I}_{d})}\left[\vec{\mu}_{c}^{\top}\vec{\hat\eta}\right]
\end{cases}
\label{eq:saddle:nohat}
\end{align}
where all the relevant quantities have been defined in Appendix \ref{subsec:appp:formal_theorem}.

\subsection{Training and generalization errors}

It now remains to compute the training and generalization errors from the free energy. Recalling that
\[\varepsilon_{\text{tr}} (\hat{\Vec{\theta}}) = \frac{1}{n}\sum\limits_{\nu=1}^{n}\ell \left(y^{\nu}, \hat{\vec{\theta}}^{\top}\vec{x}^{\nu}\right), \]
we can write using the definition of the free energy \eqref{eq:app:def_free_energy}:
\[  \lim_{n \to \infty} \varepsilon_{\text{tr}} = \lim_{\beta \to \infty} \partial_\beta f_\beta - r(\vec{\hat \theta}).\]
Computing explicitly the derivative and averaging again over the data yields
\[ \lim_{n \to \infty} \varepsilon_{\text{tr}} = -\lim_{\beta \to \infty} \partial_\beta \Psi_\ell, \]
where $\Psi_\ell$ is the free energy contribution of the loss defined in \eqref{eq:app:def_free_loss}. Writing explicitly this derivative,
\begin{equation} \label{eq:app:free_loss_derivative}
-\partial_\beta \Psi_\ell = \sum_{c \in \mathcal C} \mathbb E_{\xi_c \sim \mathcal N(0, 1)} \int_{\mathbb{R}} \dd{y} \frac{\mathcal Z_0\left(y, \omega_c^{(0)}, \rho - \frac{m_c^2}{q_c}\right)}{\mathcal Z_\ell(y, \omega_c^{(\ell)}, V_c)} \underbrace{\int_{\mathbb R} \frac{\dd{\lambda}}{\sqrt{2\pi V_c}} e^{-(\lambda - \omega_c^{(\ell)})V_c^{-1}(\lambda - \omega_c^{(\ell)}, V_c)) - \beta \ell(y, \lambda)} \ell(y, \lambda)}_{\tilde{\mathcal Z}_\ell(y, \omega_c^{(\ell)}, V_c)}
\end{equation}
where
\[ \omega_c^{(0)} = \pi_c + \frac{m_c}{\sqrt{q_c}} \xi_c \quad \text{and} \quad \omega_c^{(\ell)} = h_c + \sqrt{q_c} \xi_c. \]

We can now make the change of variables in \eqref{eq:app:beta_rescaling}, and use again Laplace's approximation: if
\[ \eta_c = \prox_{V\ell(y, \cdot)}(\omega_c^{(\ell)}) \quad \text{and} \quad \mathcal L(\lambda) = -(\lambda - \omega_c^{(\ell)})V_c^{-1}(\lambda - \omega_c^{(\ell)}) - \beta \ell(y, \lambda),\]
then
\[ \mathcal Z_\ell(y, \omega_c^{(\ell)}, V_c) \sim e^{-\beta \mathcal L(\eta_c)} \quad \text{and} \quad \tilde{\mathcal Z}_\ell(y, \omega_c^{(\ell)}, V_c) \sim e^{-\beta \mathcal L(\eta_c)} \ell(y, \eta_c),\]
and all of the overlaps will concentrate around the solutions of \eqref{eq:saddle:hat}, \eqref{eq:saddle:nohat}. Finally, the term containing $\mathcal Z_0$ is an expectation of $y$ according to $P_0(y \,\mid\, \tau_c)$,
where
\[\tau_c = \omega_c^{(0)} + \sqrt{\rho - \frac{m_c^2}{q_c}} \xi^{(0)}_c\]
and $\xi_c^{(0)}$ is a Gaussian variable independent from everything else.
Putting all together, we can write
\begin{equation}
    \lim_{n \to\infty} \varepsilon_{\text{tr}}(\vec{\hat \theta}) = \sum_{c \in \mathcal C} p_c \, \mathbb E_{\nu_c, \tau_c, y}\left[ \ell\left(y, \prox_{V_c^\star\ell(y, \cdot)}(\nu_c) \right)\right]
\end{equation} 
with
\begin{equation}
    \begin{pmatrix}
    \nu_c \\ \tau_c
    \end{pmatrix} \sim \mathcal N \left( \begin{bmatrix}
        \pi_c \\ h_c
    \end{bmatrix}, \begin{bmatrix}
        \rho_c & m_c^\star \\
        m_c^\star & q_c^\star
    \end{bmatrix}\right) \quad \text{and} \quad y \sim P_0( \cdot \,|\, \tau_c).
\end{equation}
The generalization error is much simpler to obtain: since $x_\text{new}, y_{\text{new}}$ are independent from the estimator $\vec{\hat\theta}$, we simply have
\begin{equation}
    \lim_{n \to\infty} \varepsilon_{\text{gen}}(\vec{\hat \theta}) = \sum_{c \in \mathcal C} p_c\, \mathbb E_{\nu_c, \tau_c, y}\left[ \ell\left(y, \nu_c \right)\right]
\end{equation}
where $\nu_c, \tau_c, y$ follow the same distribution as above.

\subsection{Examples}
\label{subsec:app:examples}
\paragraph{Ridge penalty} Consider a particular case of the general equations reported above: the case of a ridge penalty 
\[ r(\vec{\theta}) = \frac{\lambda}{2}||\vec{\theta}||^2_{2}. \] 
We then have:
\begin{align}
	\prox_{r}(\vec{x}) = (1+\lambda)^{-1}\vec{x}
\end{align}
\noindent which simplify the prior equations considerably. Indeed, we can now compute every expectation in \eqref{eq:saddle:nohat}, which yields the following fixed-point equations:
\begin{align}
\label{eq:saddle:l2_hat}
&\begin{cases}
	\hat{V}_{c} &= -\alpha p_{c} \mathbb{E}_{\xi_c\sim\mathcal{N}(0,1)}\left[\int\dd y~\mathcal{Z}_{0}\left(y,\pi_{c} + \frac{m_{c}}{\sqrt{q}_{c}} \xi_c,\rho-\frac{m_{c}^2}{q_{c}}\right)\partial_{\omega}f_{\ell}(y,h_{c}+\sqrt{q_{c}}\xi_c,V_{c})\right]\\
	\hat{q}_{c} &= \alpha p_{c} \mathbb{E}_{\xi_c\sim\mathcal{N}(0,1)}\left[\int\dd y~\mathcal{Z}_{0}\left(y,\pi_{c} + \frac{m_{c}}{\sqrt{q}_{c}} \xi_c,\rho-\frac{m_{c}^2}{q_{c}}\right) f_{\ell}(y,h_{c}+\sqrt{q_{c}}\xi_c,V_{c})^2\right]\\
	\hat{m}_{c} &= \alpha p_{c} \mathbb{E}_{\xi_c\sim\mathcal{N}(0,1)}\left[\int\dd y~ 	\partial_{\omega}\mathcal{Z}_{0}\left(y,\pi_{c} + \frac{m_{c}}{\sqrt{q}_{c}} \xi_c,\rho-\frac{m_{c}^2}{q_{c}}\right) f_{\ell}(y,h_{c}+\sqrt{q_{c}}\xi_c,V_{c})\right]  \\
	\hat{h}_{c} &= \alpha p_{c} \mathbb{E}_{\xi_c\sim\mathcal{N}(0,1)}\left[\int\dd y~ 	\mathcal{Z}_{0}\left(y,\pi_{c} + \frac{m_{c}}{\sqrt{q}_{c}} \xi_c,\rho-\frac{m_{c}^2}{q_{c}}\right) f_{\ell}(y,h_{c}+\sqrt{q_{c}}\xi_c,V_{c})\right] 
\end{cases} \\[1em]
\label{eq:saddle:l2_nohat}
&\begin{cases}
V_c &=  \frac{1}{d}\tr\left[\Sigma_c \left(\lambda I_d + \hat\Sigma\right)^{-1}\right]\\ \noalign{\vskip 0.2em}
q_c &=\frac{1}{d}\tr\left[\left(\sum\limits_{c'\in \mathcal{C}}\hat{q}_{c'}\Sigma_{c'}+\sum\limits_{c',c'' \in \mathcal{C}} \vec{\hat \mu}_{c'}\vec{\hat \mu}_{c''}^\top  \right)\Sigma_c\left(\lambda I_d + \hat\Sigma\right)^{-2}\right]\\ \noalign{\vskip 0.2em}
m_c &= \frac{1}{d}\tr \left[ \left(\sum\limits_{c'\in \mathcal{C}}\vec{\hat \mu}_{c'}\vec{\theta}_0^\top \right) \Sigma_c \left(\lambda I_d + \hat\Sigma\right)^{-1}\right]\\ \noalign{\vskip 0.2em}
h_c &=\frac{1}{d}\tr\left[\left(\sum\limits_{c'\in \mathcal{C}}\vec{\hat \mu}_{c'}\vec{\mu}_c^\top \right)\left(\lambda I_d + \hat\Sigma\right)^{-1}
\right] 
\end{cases} 
\end{align}

\paragraph{Gaussian covariate model}

The equations for the Gaussian covariate model can be found in \cite{gcm_paper}; they also correspond to taking $|\mathcal C| = 1$ in the ones above. We reproduce them here for completeness:
\begin{align}
\label{eq:saddle:gcm_l2_hat}
&\begin{cases}
	\hat{V} &= -\alpha\mathbb{E}_{\xi\sim\mathcal{N}(0,1)}\left[\int\dd y~\mathcal{Z}_{0}\left(y,\pi + \frac{m}{\sqrt{q}}\xi,\rho-\frac{m^2}{q}\right)\partial_{\omega}f_{\ell}(y,h+\sqrt{q}\xi,V)\right]\\
	\hat{q} &= \alpha\mathbb{E}_{\xi\sim\mathcal{N}(0,1)}\left[\int\dd y~\mathcal{Z}_{0}\left(y,\pi + \frac{m}{\sqrt{q}}\xi,\rho-\frac{m^2}{q}\right) f_{\ell}(y,h+\sqrt{q}\xi,V)^2\right]\\
	\hat{m} &= \alpha\mathbb{E}_{\xi\sim\mathcal{N}(0,1)}\left[\int\dd y~ 	\partial_{\omega}\mathcal{Z}_{0}\left(y,\pi + \frac{m}{\sqrt{q}}\xi,\rho-\frac{m^2}{q}\right) f_{\ell}(y,h+\sqrt{q}\xi,V)\right]  \\
	\hat{h} &= \alpha\mathbb{E}_{\xi\sim\mathcal{N}(0,1)}\left[\int\dd y~ 	\mathcal{Z}_{0}\left(y,\pi + \frac{m}{\sqrt{q}}\xi,\rho-\frac{m^2}{q}\right) f_{\ell}(y,h+\sqrt{q}\xi,V)\right] 
\end{cases} \\[1em]
\label{eq:saddle:gcm_l2_nohat}
&\begin{cases}
V &=  \frac{1}{d}\tr\left[\Sigma \left(\lambda I_d + \hat\Sigma\right)^{-1}\right]\\ \noalign{\vskip 0.2em}
q &=\frac{1}{d}\tr\left[\left(\hat{q}\Sigma+ \vec{\hat \mu}\vec{\hat \mu}^\top  \right)\Sigma\left(\lambda I_d + \hat\Sigma\right)^{-2}\right]\\ \noalign{\vskip 0.2em}
m &= \frac{1}{d}\tr \left[ \vec{\hat \mu}\vec{\theta}_0^\top \Sigma \left(\lambda I_d + \hat\Sigma\right)^{-1}\right]\\ \noalign{\vskip 0.2em}
h &=\frac{1}{d}\tr\left[\vec{\hat \mu}\vec{\mu}^\top\left(\lambda I_d + \hat\Sigma\right)^{-1}
\right] 
\end{cases} 
\end{align}
where this time
\begin{equation}
    \vec{\hat \mu} = \hat h \vec{\mu} + \hat m \Sigma \vec{\theta}_0 \qquad \text{and} \qquad \hat \Sigma = \hat V \Sigma.
\end{equation}
The errors are given by
\begin{align}
 \lim_{n \to\infty} \varepsilon_{\text{tr}}(\vec{\hat \theta}) &= \mathbb E_{\nu, \tau, y}\left[ \ell\left(y, \prox_{V^\star\ell(y, \cdot)}(\nu) \right)\right] \\
    \lim_{n \to\infty} \varepsilon_{\text{gen}}(\vec{\hat \theta}) &= \mathbb E_{\nu, \tau, y}\left[ \ell\left(y, \nu \right)\right]
\end{align}
with
\begin{equation}
    \begin{pmatrix}
    \nu \\ \tau
    \end{pmatrix} \sim \mathcal N \left( \begin{bmatrix}
        \pi \\ h
    \end{bmatrix}, \begin{bmatrix}
        \rho & m^\star \\
        m^\star & q^\star
    \end{bmatrix}\right) \quad \text{and} \quad y \sim P_0( \cdot \,|\, \tau).
\end{equation}

\paragraph{Ridge regression} We place ourselves in the ridge regression case, where
\begin{align}
    P_0(\tau|\Delta) &= \mathcal{N}(\tau,\Delta) & l(y,\hat{y}) &= (y-\hat{y})^2.
\end{align}
In this case, the equations in \eqref{eq:saddle:l2_nohat} simplify even further, yielding
\begin{align}
\label{eq:saddle:ridge_hat}
&\begin{cases}
	\hat{V}_{c} &= \frac{\alpha p_{c}}{1+V_c}\\
	\hat{q}_{c} &= \frac{\alpha p_c}{(1+V_c)^2}(\rho_c + \Delta + q_c -2m_c +(h_c-\pi_c)^2)\\
	\hat{m}_{c} &= \frac{\alpha p_{c}}{1+V_c}  \\
	\hat{h}_{c} &= \frac{\alpha p_c (\pi_c - h_c)}{1+V_c}
\end{cases} \\[1em]
\label{eq:saddle:ridge_nohat}
&\begin{cases}
V &=  \frac{1}{d}\tr\left[\Sigma \left(\lambda I_d + \hat\Sigma\right)^{-1}\right]\\ \noalign{\vskip 0.2em}
q &=\frac{1}{d}\tr\left[\left(\hat{q}\Sigma+ \vec{\hat \mu}\vec{\hat \mu}^\top  \right)\Sigma\left(\lambda I_d + \hat\Sigma\right)^{-2}\right]\\ \noalign{\vskip 0.2em}
m &= \frac{1}{d}\tr \left[ \vec{\hat \mu}\vec{\theta}_0^\top \Sigma \left(\lambda I_d + \hat\Sigma\right)^{-1}\right]\\ \noalign{\vskip 0.2em}
h &=\frac{1}{d}\tr\left[\vec{\hat \mu}\vec{\mu}^\top\left(\lambda I_d + \hat\Sigma\right)^{-1}
\right] 
\end{cases} 
\end{align}
The training and generalization error also benefit from a very simple expression
\begin{align}\label{eq:app:ridge_errors}
    \varepsilon_{\text{tr}} = \sum_c \frac{\hat{q}_c^\star}{\alpha} && \varepsilon_{\text{tr}} = \sum_c \frac{\hat{q}^\star_c (1+V_c^\star)^2}{\alpha}.
\end{align}
\section{Universality of Gaussian Mixture Models}
\label{sec:appendix:proofs}
In this section we prove the main results shown in  in Sec.~\ref{sec:main:theory} regarding universality properties of Gaussian Mixture Models. 
\subsection{Mean Universality : proof of Proposition \ref{prop:mean_univers}}
Let us rewrite the assumptions here. Let $\{\hat{V}^{\star}_{c}\}_{c=1}^K$ be the fixed points of the (replica) saddle point equations describing the centered Gaussian mixture problem, i.e. the solutions of of eqs.~\eqref{eq:saddle:l2_hat}, \eqref{eq:saddle:l2_nohat}.  We assume that the teacher vector and the data structure respect the following:
 \begin{align}
     &\lim_{n, d \to \infty} \frac{\Vec{\theta}_0^{\top}\Vec{\mu}_c}{d} = 0 \qquad \forall c \in \mathcal{C} \\
     &\lim_{n, d \to \infty} \frac1{d} \Vec{\theta}_0 ^{\top}\Sigma_{c'}\left(\lambda +  \sum_{c \in \mathcal C}\hat{V}^\star_{c}\Sigma_c\right)^{-1}\vec{\mu}_c \to 0 \qquad \forall (c,c') \in \mathcal{C}\times\mathcal{C} 
 \end{align}
and the loss and teacher are both symmetric:
\begin{align}
    \ell(x,y) &= \ell(-x,-y) \\
    P_0\left(y|\tau\right) &=  P_0\left(-y|-\tau\right).
\end{align}
In the saddle-point equations \eqref{eq:saddle:l2_hat}, \eqref{eq:saddle:l2_nohat}, the mean vectors appear always coupled with the overlaps $\{h_c,\hat{h}_c\}_{c \in \mathcal{C}}$. Hence, to prove Prop.~\ref{prop:mean_univers} it suffices to prove the following result:
\begin{lemma}
If $(a + b + c)$ hold, then $\{h_c = \hat{h}_c = 0\}_{c \in \mathcal{C}}$ is a fixed point for the problem.
\end{lemma}
First consider the update equations for $\{h_c\}_{c \in \mathcal{C}}$:
\begin{align}
    h_c &=\frac{1}{d}\sum_{l}\text{tr}\left[ \left(\hat{h}_l\mu_l\mu_c^\top+\hat{m}_l\Sigma_l\vec{\theta}_0\mu_c^\top\right) \left(\lambda\mat{I}_{d}+\sum_{l'}\hat{V}_{l'}\Sigma_{l'}\right)^{-1}\right] \\
    &\underset{a+c}{\simeq} \quad\frac{1}{d}\sum_{l}\text{tr}\left[ \left(\hat{h}_l\mu_l\mu_c^\top\right) \left(\lambda\mat{I}_{d}+\sum_{l'}\hat{V}_{l'}\Sigma_{l'}\right)^{-1} \right]
    \label{eq:h_eq} 
\end{align}
By continuity of the saddle-point equations, if at the fixed point $\{\hat{h}_c^\star \to 0\}_{c \in \mathcal{C}}$ holds, we also easily have that $\{h_c^\star \to 0\}_{c \in \mathcal{C}}$. 
Now assume $\{h_c^\star = 0\}_{c \in \mathcal{C}}$ at the fixed point. We exploit symmetry argument inherent to the update functions to show that $\{\hat{h}_c^\star =0\}_{c\in \mathcal{C}}$ under weak assumptions on the teacher and loss functions. We write the updates using assumption c):
\begin{align}
    \hat{h}_c &= \alpha  \mathbb{E}_{\xi\sim\mathcal{N}(0,1)}\left[\int\dd y~ 	\mathcal{Z}_{0}\left(y,\frac{m_c}{\sqrt{q_c}}\xi+\pi_c,\rho_c-\frac{m_c^2}{q_c}\right) f_{\ell}(y,h_c+\sqrt{q_c}\xi,V)\right] \\
    &\underset{h_c \to 0}{\simeq} \quad \alpha  \mathbb{E}_{\xi\sim\mathcal{N}(0,1)}\left[\int\dd y~ 	\mathcal{Z}_{0}\left(y,\frac{m_c}{\sqrt{q_c}}\xi,\rho_c-\frac{m_c^2}{q_c}\right) f_{\ell}(y,\sqrt{q_c}\xi,V_c)\right] \label{eq:app:hat_hc_odd}
\end{align}

Reminding the definition of the teacher measure term: 
\begin{align}
    Z_0(y,\omega,V) = \mathbb{E}_{\lambda \sim \mathcal{N}(\omega,V)} \left[P_0(y|\lambda)\right] && \lambda \equiv \frac{\Vec{x}^\top\Vec{\theta}_0}{\sqrt{d}}
\end{align}
\noindent we see that the symmetry conditions impose basic restrictions on the label generation:
\begin{align}
     f) \quad \mathbb{E}_{\lambda \sim \mathcal{N}(\omega,V)} \left[P_0(y|\lambda)\right] =  \mathbb{E}_{\lambda \sim \mathcal{N}(-\omega,V)} \left[P_0(-y|\lambda)\right]
\end{align}
and hence $\mathcal Z_0$ is even in its second argument. Additionally, we have
\begin{align}
    f_\ell(y,\omega,V) = \frac{1}{V}(\prox\,_{Vl(y,\cdot)}(\omega) - \omega ) && \prox\,_{V\ell(y,\cdot)} (w) = \argmin_{z} \left[\frac{1}{2V}(z-w)^2 + \ell(y,z)\right]
\end{align}
The symmetry condition on $\ell$ then implies that
\begin{align}
    \prox\,_{V\ell(-y,\cdot)}(-\omega) = - \prox\,_{V\ell(y,\cdot)}(\omega),
\end{align}
so $f_\ell$ is odd in its second argument. All that's left to notice is that \eqref{eq:app:hat_hc_odd} is an Gaussian integral of an odd function, hence it is equal to 0.

\subsection{Gaussian Universality : proof of Proposition \ref{prop:gauss_univ}}
\label{sec:app:gauss_univ_proof}
We proved that the fixed point respects $\{h_c^\star,\hat{h}_c^\star= 0,0\}_{c\in \mathcal{C}}$ under the hypothesis of Prop.~\ref{prop:mean_univers}.
Assume now that the mixture is homogeneous:
\begin{align}
    \Sigma_c = \Sigma \qquad \forall c \in \mathcal{C}
\end{align}
Then we have
\begin{align}
    \varepsilon^{\rm GMM}_{gen} \left(\{\Vec{\mu}_c\}_{c=1}^K,\{\Sigma\}_{c=1}^K\right) &\simeq \varepsilon^{\rm GMM}_{gen}\left(\Vec{0},\{\Sigma\}_{c=1}^K\right) \\
    \varepsilon^{\rm GMM}_{tr} \left(\{\Vec{\mu}_c\}_{c=1}^K,\{\Sigma\}_{c=1}^K\right) &\simeq \varepsilon^{\rm GMM}_{tr}\left(\Vec{0},\{\Sigma\}_{c=1}^K\right)
\end{align}
But the right-hand side of those equations corresponds to a distribution of the form
\[\vec{x} \sim \sum_{c\in\mathcal C} \mathcal N(\vec{0}, \Sigma) = \mathcal N(\vec{0}, \Sigma) !\]
This proves Proposition \ref{prop:gauss_univ}

\subsection{Covariance universality} Surprisingly in the limit of vanishing regularization, $\lambda \to 0^+$, the covariance is not relevant for the high dimensional learning problem. We show that when a unique minimizer of the loss exists, and we can take safely the limit $\lambda \to 0^+$ in the saddle point equations, the covariance $\Sigma$ disappears completely from the overlap expression. Indeed assuming the minimizer $\vec{\hat \theta}$ is unique, we can safely simplify expression of the type: 
\begin{align}
\label{eq:app:no_lamb_simplific}
    \lim_{\lambda \to 0^+} \Sigma^n(\lambda \mat{I}_d + \hat{V}\Sigma)^{n^\prime} = \frac{1}{\hat{V}}\Sigma^{n-n^\prime}
\end{align}
Then by plugging in the $\lambda \to 0^+$ simplification in eqs.~\eqref{eq:saddle:gcm_l2_hat},
\eqref{eq:saddle:gcm_l2_nohat} we have complete independence from the covariance matrix $\Sigma$:
\begin{align}
&\begin{cases}
\hat{V} &= -\alpha \mathbb{E}_{\xi\sim\mathcal{N}(0,1)}\left[\int\dd y~\mathcal{Z}_{0}\left(y,\frac{m_{c}}{\sqrt{q_{c}}}\xi,\rho_{c}-\frac{m_{c}^2}{q_{c}}\right)\partial_{\omega}f_{\ell}(y,\sqrt{q_{c}}\xi,V_{c})\right]\\
\hat{q} &= \alpha \mathbb{E}_{\xi\sim\mathcal{N}(0,1)}\left[\int\dd y~\mathcal{Z}_{0}\left(y,\frac{m_{c}}{\sqrt{q}}\xi,\rho-\frac{m^2}{q}\right) f_{\ell}(y,\sqrt{q}\xi,V)^2\right]\\
\hat{m} &= \alpha \mathbb{E}_{\xi\sim\mathcal{N}(0,1)}\left[\int\dd y~ 	\partial_{\omega}\mathcal{Z}_{0}\left(y,\frac{m}{\sqrt{q}}\xi,\rho-\frac{m^2}{q}\right) f_{\ell}(y,\sqrt{q}\xi,V)\right] 
\end{cases} \\
&\begin{cases}
    V &= 1 \\
    q &= \hat{q} + \rho \hat{m}^2 \\
    m &= \hat{m} \rho
\end{cases}
\end{align} 
which concludes the proof.
\subsection{Strong universality}
We never made any specific assumption on the loss up to now, apart the very general symmetry condition in Assumption \ref{assump:symmetric}. In this section we show strong universality for square loss regression for \textit{any GMM} estimation problem respecting the mean universality property (Prop.~\ref{prop:mean_univers}). 
We assume to consider a ridge regression problem in the underparametrized regime $\alpha > 1$:
\begin{align}
    P_0(\tau|\Delta) = \mathcal{N}(\tau,\Delta) && \ell(y,\hat{y}) = (y-\hat{y})^2 
\end{align}
the replicas then correspond to equations \eqref{eq:saddle:ridge_hat}, \eqref{eq:saddle:ridge_nohat}. In the $\lambda \to 0^+$ limit, the equation simplify greatly: 
\begin{align}
\label{eq:app:strong_univ_sp_hats}
&\begin{cases}
	\hat{V}_{c} &= \frac{\alpha p_{c}}{1+V_c} \\
	\hat{q}_{c} &= \frac{\alpha p_c}{(1+V_c)^2}(\rho_c + \Delta + q_c -2m_c) \\
 \hat{m}_{c} &= \frac{\alpha p_{c}}{1+V_c}
\end{cases} \\
\label{eq:strong_univ_sp}
&\begin{cases}
V_c &=  \frac{1}{d}\tr\left[\Sigma_c (\sum_{l \in \mathcal{C}} \hat{V}_l \Sigma_l)^{-1}\right]\\
q_c &= \frac{1}{d}\sum_{l \in \mathcal{C}}^K\tr\left[\hat{q}_l\Sigma_l\Sigma_c (\sum_{l^{\prime} \in \mathcal{C}} \hat{V}_{l^{\prime}} \Sigma_{l^{\prime}})^{-2}\right]  +  \frac{1}{d}\sum_{(l,l)'\in \mathcal{C} \times \mathcal{C}}\tr\left[\hat{m}_l\hat{m}_{l'}\Sigma_l\vec{\theta}_0\Vec{\theta}_0^\top\Sigma_{l'}\Sigma_c (\sum_{l^{\prime \prime} \in \mathcal{C}} \hat{V}_{l^{\prime \prime}} \Sigma_{l^{\prime \prime}})^{-2}\right] \\
m_c &= \frac{1}{d}\sum_{l=1}^K\tr \left[ \hat{m}_l\Sigma_l\vec{\theta}_0\vec{\theta}_0^T \Sigma_c (\sum_{l'} \hat{V}_{l'} \Sigma_{l'})^{-1}\right]
\end{cases} 
\end{align}

Consider the equation for the overlaps $\{q_c\}_{c \in \mathcal{C}}$:
\begin{align}
    q_c &= \frac{1}{d}\sum_{l=1}^K\tr\left[\hat{q}_l\Sigma_l\Sigma_c  (\sum_{l'} \hat{V}_{l'} \Sigma_{l'})^{-2}\right]  +  \frac{1}{d}\sum_{l,k=1}^K \tr\left[\hat{m}_l \hat{m}_k\Sigma_l \Vec{\theta}_0\vec{\theta}_{0}^T\Sigma_k\Sigma_c (\sum_{l'} \hat{V}_{l'} \Sigma_{l'})^{-2} \right] \\
    & \equiv R_c + G_c
\end{align}
The error metrics can be decomposed as:
\begin{align}
\label{eq:app:et_strong_univ}
    \varepsilon_{\text{tr}} &= \sum_{c\in\mathcal{C}}\frac{ p_c}{(1+V_c)^2}(R_c+\Delta) + \sum_{c \in \mathcal{C}}\frac{ p_c}{(1+V_c)^2}(\rho_c + G_c -2m_c)  \\
    \label{eq:app:eg_strong_univ}
    \varepsilon_{\text{gen}} &= \sum_{c\in\mathcal{C}} p_c(R_c+\Delta) + \sum_{c \in \mathcal{C}} p_c(\rho_c + G_c -2m_c ) 
\end{align}

We focus on the second term and show that it is equal to zero. 
Looking at eqs.~\eqref{eq:strong_univ_sp},\eqref{eq:app:strong_univ_sp_hats}, we note that:
\begin{align}
    \hat{m}_c = \hat{V}_c \qquad \forall c \in \mathcal{C}
    \label{eq:app:mhat=vhat}
\end{align}
By using eq.~\eqref{eq:app:mhat=vhat} we can simplify the equations for $\{\hat{m}_c,G_c\}$:
\begin{align}
    G_c = \rho_c, \qquad m_c = \rho_c \qquad \forall c \in \mathcal{C}
\end{align}
Plugging this relations in the saddle point equations we obtain:
\begin{align}
&\begin{cases}
	\hat{V}_{c} &= \frac{\alpha p_{c}}{1+V_c} = \hat{m}_{c}\\
	\hat{q}_{c} &= \frac{\alpha p_c}{(1+V_c)^2}(\Delta + R_c) 
\end{cases} \\
&\begin{cases}
V_c &=  \frac{1}{d}\tr\left[\Sigma_c (\sum_{l \in \mathcal{C}} \hat{V}_l \Sigma_l)^{-1}\right]\\
q_c &= \frac{1}{d}\sum_{l=1}^K\tr\left[\hat{q}_l\Sigma_l\Sigma_c  (\sum_{l'} \hat{V}_{l'} \Sigma_{l'})^{-2}\right] + \rho_c \equiv R_c +  \rho_c \\
m_c &= \rho_c
\end{cases} 
\end{align}
The fixed point of the equations does not depend on $\{\rho_c\}_{\in \mathcal{C}}$ anymore and we are left with equations only for $\{\hat{V}_c,\hat{q}_c,V_c,R_c\}_{c \in \mathcal{C}}$: 
\begin{align}
&\begin{cases}
	\hat{V}_{c} &= \frac{\alpha p_{c}}{1+V_c}\\
	\hat{q}_{c} &= \frac{\alpha p_c}{(1+V_c)^2}(\Delta + R_c) 
\end{cases} \\
&\begin{cases}
V_c &=  \frac{1}{d}\tr\left[\Sigma_c (\sum_{l \in \mathcal{C}} \hat{V}_l \Sigma_l)^{-1}\right]\\
R_c &= \frac{1}{d}\sum_{l=1}^K\tr\left[\hat{q}_l\Sigma_l\Sigma_c  (\sum_{l'} \hat{V}_{l'} \Sigma_{l'})^{-2}\right] 
\end{cases} 
\end{align}
The errors can be computed in a closed form with a series of algebraic manipulation.
One can verify by algebraic manipulation the following relations:
\begin{align}
\label{eq:app:1st_relation_strong_univ}
    \sum_c V_c \hat{V}_c &= \sum_c\frac{1}{d}\tr\left[\hat{V}_c\Sigma_c (\sum_{l \in \mathcal{C}} \hat{V}_l \Sigma_l)^{-1}\right] = 1 \\
    \sum_cR_c\hat{V}_c - V_c \hat{q}_c &=   \frac{1}{d}\sum_{l=1}^K\tr\left[\hat{q}_l\Sigma_l (\sum_c\hat{V}_c\Sigma_c)  (\sum_{l'} \hat{V}_{l'} \Sigma_{l'})^{-2}\right]   - \frac{1}{d}\tr\left[\sum_c\hat{q}_c\Sigma_c (\sum_{l \in \mathcal{C}} \hat{V}_l \Sigma_l)^{-1}\right] = 0
\end{align}
Now we express everything in eq.~\eqref{eq:app:1st_relation_strong_univ} in terms of the non-hatted overlaps to obtain:
    \begin{align}
         &\star)\, \sum_c \frac{\alpha p_c}{1+V_c}V_c = 1
         \\
        &\#) \,\sum_c \frac{\alpha p_c}{(1+V_c)^2} R_c  =  
        \sum_c \frac{\alpha p_c}{(1+V_c)^2} V_c \Delta 
    \end{align}
We remark that we can write the training and generalization error can be written as:
\begin{align}
\label{eq:app:metrics_definition}
    \varepsilon_{\text{tr}} = \sum_c \frac{\hat{q}_c}{\alpha} && \varepsilon_{\text{tr}} = \sum_c \frac{\hat{q}^\star_c (1+V_c^\star)^2}{\alpha}
\end{align}
So if we compute the quantity $\sum_c \hat{q}_c$ we conclude. By plugging in $(\#)$ into the expression above we obtain:
    \begin{align}
        \sum_c\hat{q}_c &= \sum_c \frac{\alpha p_c}{(1+V_c)^2}(R_c+\Delta) \\
        &\underset{(\#)}{=} \Delta \sum_c \frac{\alpha p_c}{(1+V_c)^2}+  \Delta \sum_c \frac{\alpha p_c}{(1+V_c)^2}V_c   \\
        &= \Delta \sum_c \frac{\alpha p_c}{1+V_c} 
    \end{align}
and finally using relation $(\star)$ we prove the theorem:
    \begin{align}
        \varepsilon_{\text{tr}} &= \sum_c \frac{\hat{q}_c^\star}{\alpha} \\
        &=  \sum_c \frac{\hat{q}_c^\star}{\alpha} \pm \frac{\Delta}{\alpha} 
        \\
        &\underset{(\star)}{=} \Delta (1-\frac{1}{\alpha}) 
    \end{align}
On the other hand, the generalization error does not respect the strong universlaity statement as we analyze in the next section (See Fig.~\ref{fig:hetero_mixture}).

\section{Non-universality of Gaussian Mixture Models}
\label{sec:appendix:correlated}
In the previous section we enumerated a series of results unveiling universality of GMM. Now we want to study the dual task: when GMM model \textit{are not} universal?  
We have two ways to break universality: a) allow strong heterogeneity in the data structure; b) consider labels which are strongly correlated with the data structure. The plan for this section is to review more in detail these two processes for universality breaking. 

\subsection{Strongly heterogeneous mixtures}
\begin{figure*}
\centering
\includegraphics[width=0.8\textwidth]{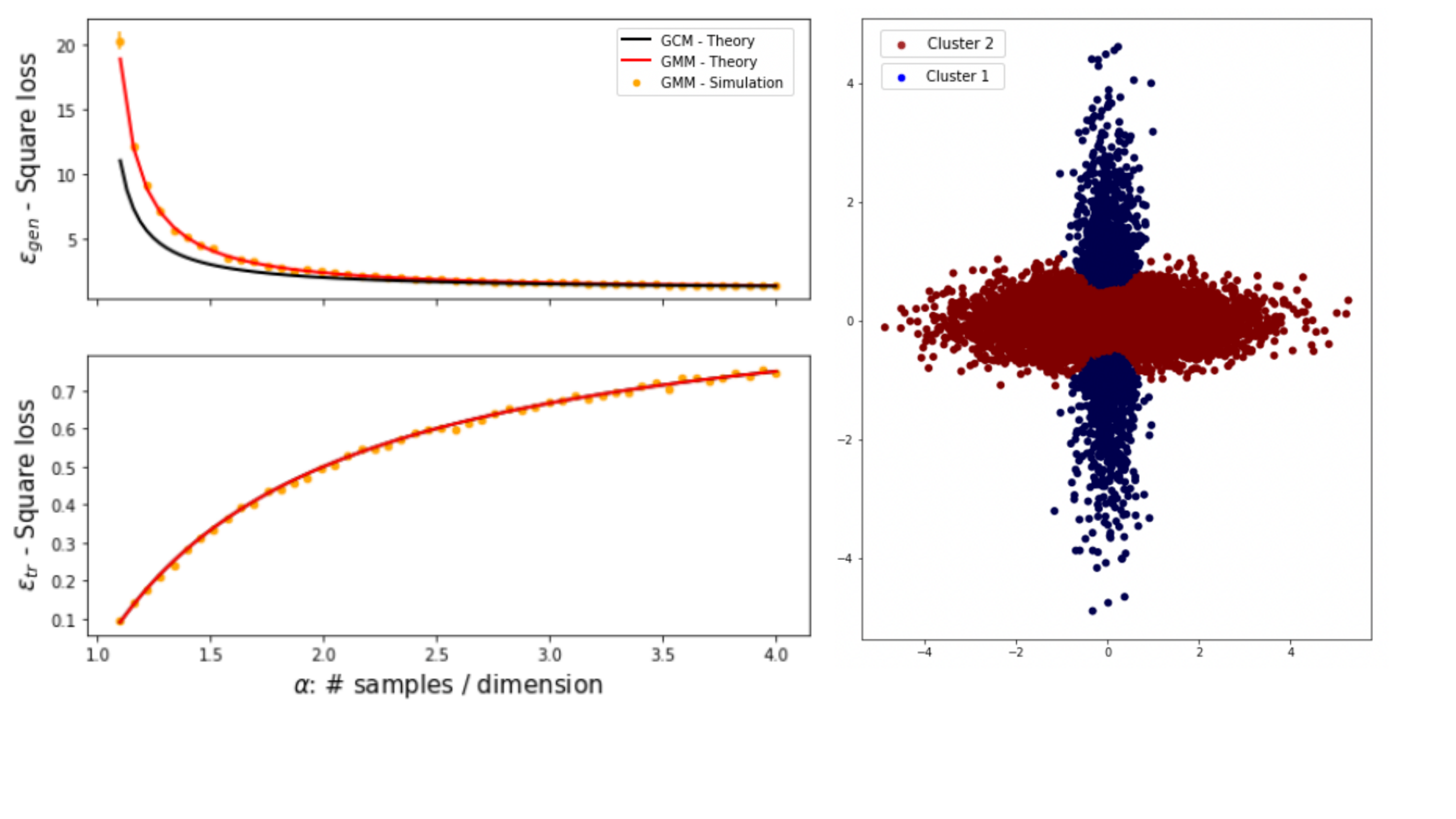}
\vspace{-3em}
\caption{Training and generalization error for ridge regression on the GMM defined in eq.~\eqref{eq:hetero_mixt}. \textbf{Left}: We compare the performance of the Gaussian asymtptotics (solid black line) with the GMM one (solid red line), while the orange dots represents simulations which agrees with the theoretical predictionas predicted by Theorem.~\ref{prop:exact_asymp_gmm}. \textbf{Right}: Two dimensional toy plot of $10000$ samples coming from the heterogeneous Gaussian mixture.  }
\label{fig:hetero_mixture}
\end{figure*}

We proved in Theorem.~\ref{prop:ridge_universality} a strong universality statement for the training loss of ridge regression. However in this section we want to clarify that the theorem is not valid beyond its assumption for general mixtures. We prove this by analyzing a counterexample: a strongly heterogeneous 2-clusters Gaussian Mixture:
\begin{align}
\label{eq:hetero_mixt}
    \Sigma_+ &= \text{diag}(0.1,\dots,0.1,1.9,\dots,1.9) \qquad p_+ = 0.8 \\
    \Sigma_- &= \text{diag}(1.9,\dots,1.9,0.1,\dots,0.1) \qquad p_- = 0.2
\end{align}
We present in Fig.~\ref{fig:hetero_mixture} the comparison of the heterogeneous GMM performance with the Gaussian theory: although the training errors coincide as predicted by Theorem.~\ref{prop:ridge_universality}, the generalization errors are different. However, we remark that real data after preprocessing seem homogeneous enough to obtain a good agreement with the exact Gaussian asymtptotics as we see in Fig.~\ref{fig:fig1}.
\subsection{Correlated teachers}
In this section we investigate in deeper detail the controlled setting of Theorem.~\ref{prop:non_univ_mixt}. We  do not necessarily consider estimation problems which respects the assumption of mean universality property (Prop.~\ref{prop:mean_univers}), and in fact we precisely analyze the consequences if we relax these conditions. We perform ridge regression on a two-mixtures GMM with opposite means and same covariance matrix proportional to $\mat{I}_d$. 
In this scenario the exact asymptotics for the performance of the ERM estimator admits a closed form expression.
We remind here the main parameters for the theoretical analysis are:
\begin{align}
    \rho = \frac{1}{d}\Vec{\theta}_0^\top\Sigma\Vec{\theta}_0= \frac{1}{d}\lVert\Vec{\theta}_0\rVert_2^2 \qquad 
    \gamma = \frac{1}{d} ||\Vec{\mu}||_2^2 \qquad 
    \Delta = \mathbb{E}_\xi \xi^2  \qquad
    \vec{\hat \mu_{\pm}} = \hat h_{\pm} \vec{\mu}_{\pm} + \hat m_{\pm} \Sigma \vec{\theta}_0
\end{align}
In this simple setting we can simplify the general mixture equations in eqs.~\eqref{eq:saddle:gcm_l2_hat},\eqref{eq:saddle:gcm_l2_nohat} without doing any assumption on the teacher. We first map the equation to a single Gaussian problem $\mathcal{G}_{\Vec{\theta},\mat{I}_d,\Vec{\mu}}$. We have to slightly modify the proof of Prop.~\ref{prop:gauss_univ} presented in Sec.~\ref{sec:app:gauss_univ_proof}, indeed we cannot blindly assume that the mean overlaps $\{h_+,h_-\}$ are zero but one can show that at the fixed point the overlaps respect:
\begin{align}
    \hat{h}_+ ^\star= - \hat{h}_-^\star, \qquad \hat{m}_+^\star = \hat{m}_-^\star \qquad \to  \qquad \Vec{\hat\mu}_+^\star = \Vec{\hat\mu}_-^\star
\end{align}
\begin{align}
    \hat{V} = \hat{V}_+ +  \hat{V}_- && \hat{q} = \hat{q}_+ + \hat{q}_- && \hat{m} = \hat{m}_+ + \hat{m}_- && \hat{h} = \hat{h}_+ - \hat{h}_-   \\ 
    V = V_+ = V_- && q = q_+ = q_- && m = m_+ = m_- && h = h_+ = - h_-
\end{align}
the fixed point of the replica equations for the mixture defined above are mapped to the one of a single Gaussian problem for the overlaps $\{V,m,q,h,\hat{V},\hat{m}, \hat{q} , \hat{h}  \}$.
Now we can simplify them even further by plugging in the assumption on the covariance, all the traces simplify and we get: 
\begin{equation}\label{eq:sigma_simplify}
    \begin{cases}
		V &= \frac1{\lambda +  \hat V}\\
		q &= V^2\left( \hat q +  \rho \hat m^2 + 2\pi  \hat m \hat h + \gamma \hat h^2 \right)\\
		m &= V\left( \rho \hat m +  \pi \hat h\right) \\
		h &= V\left(\gamma \hat h + \pi  \hat m\right)
	\end{cases}
\end{equation}
These equations are actually solvable! Define
\[\eta := V\hat V = \frac{\alpha V}{1+V};\]
Now, we know that the generalization error satisfies:
\[ V^2 \hat q = \frac{V^2\alpha}{(1+V)^2}\varepsilon_{\text{gen}} = \frac{\eta^2}{\alpha} (\varepsilon_{\text{gen}} ). \]
We can plug this into the equation for $q$ to get
\begin{align}
    \varepsilon_{\text{gen}} &= \rho + \underbrace{\left(\frac{\eta^2}{\alpha} (\varepsilon_{\text{gen}} ) + \rho (V\hat m)^2 + 2\pi (V \hat m)(V\hat h) + \gamma (V\hat h)^2\right)}_{q} - 2\underbrace{\left( \rho V\hat m + \pi V \hat h\right)}_m + (\pi - h)^2 \\
    \varepsilon_{\text{tr}} &= \frac{\varepsilon_{\text{gen}}}{(1+V)^2}
\end{align}
with the relations:
\[ V = \frac{\eta}{\alpha-\eta},\qquad V\hat m = \eta, \quad V \hat h = \eta  (\pi - h)  \]
and the expression of $\pi - h$ as a function of $\eta$:
\begin{align}
    \pi - h = \pi - V(\gamma \hat{h} + \pi \hat{m}) = \pi - \gamma \eta (\pi - h) - \pi \eta \iff \pi - h = \pi\frac{1-\eta}{1+\gamma\eta}
\end{align}
we simplify everything in terms only of $\eta$ to get:
\begin{align}
\label{eq:gen_error_full generality}
    \varepsilon_{\text{gen}} &= \frac{\alpha\left(\Delta +(\eta -1)^2 \rho \right)}{\left(\alpha  - \eta ^2\right)}  \left( 1 - \pi^2 \frac{(\eta -1)^2 \left(\gamma  \eta ^2+2 \eta -1\right)}{\left(\Delta +(\eta -1)^2 \rho \right)(1+\gamma \eta)^2}\right)  \\ 
 \label{eq:train_error_full generality}   
    \varepsilon_{\text{tr}} &= \frac{(\alpha-\eta)^2\left(\Delta +(\eta -1)^2 \rho \right)}{\alpha\left(\alpha - \eta^2\right) }  \left(1 - \pi^2 \frac{(\eta -1)^2 \left(\gamma  \eta ^2+2 \eta -1\right)}{\left(\Delta +(\eta -1)^2 \rho \right)(1+\gamma \eta)^2}\right)
\end{align}
All that remains is to solve for $\eta$ using the equations for $V$ and $\hat V$, which yields
\begin{equation}\label{eq:eta_expression}
\eta = 1 + \frac12 \left( \alpha - 1 + \lambda - \sqrt{4\lambda + (\alpha -1 + \lambda)^2} \right). 
\end{equation}
 
\paragraph{Vanishing regularization}
We can take $\lambda = 0$ inside \eqref{eq:eta_expression} even when $\alpha < 1$, and we find
\begin{align}
\label{eq:no_lamb_eta}
    \eta(\lambda = 0) =  \min(\alpha, 1)
\end{align}
Finally, this yields
\begin{align}
   \varepsilon_{gen} &= \begin{cases}
 \frac{\alpha  \Delta }{\alpha -1}  \qquad &\alpha \geq 1 \\
 \frac{\pi^2(\alpha -1) \left(\alpha ^2 \gamma +2 \alpha -1\right)}{(\alpha  \gamma +1)^2}-\frac{(\alpha -1)^2 \rho +\Delta }{\alpha -1} \qquad &\text{else} 
\end{cases}  \\
 \varepsilon_{\text{train}} &= 
 \begin{cases}
 \frac{(\alpha -1)  \Delta }{\alpha} \qquad &\alpha \geq 1 \\
 0 \qquad &\text{else} \\
\end{cases}
\end{align}
We retrieve the results of Corollary.~\ref{prop:corr_teach_univers}: even with correlated teachers we show that in the underparametrized regime we have Gaussian universality. 

\paragraph{Extension of Corollary.~\ref{prop:corr_teach_univers}:} The Gaussian universality result for correlated teacher at vanishing regularization can be extended as well to general balanced mixtures with homoscedastic covariance, i.e. the case where
\[ p_c = \frac1{|\mathcal C|}, \quad \Sigma_c = \Sigma, \quad \lambda = 0. \]

In order to do so, we consider the general saddle point equation for ridge regression in eq.~\eqref{eq:saddle:ridge_hat},\eqref{eq:saddle:ridge_nohat}, and plug in the assumptions:
\begin{align}
&\begin{cases}
	\hat{V} &= \frac{\alpha}{|\mathcal C|(1+V)}\\
	\hat{q}_{c} &= \frac{\alpha}{|\mathcal C|(1+V)^2}(\rho + \Delta + q_c -2m +(h_c-\pi_c)^2)\\
	\hat{m} &= \frac{\alpha}{|\mathcal C|(1+V)}  \\
	\hat{h}_{c} &= \frac{\alpha (\pi_c - h_c)}{|\mathcal C|(1+V)}
\end{cases} \\
&\begin{cases}
V &=  \frac{1}{d}\tr\left[\Sigma \hat\Sigma^{-1}\right]\\ \noalign{\vskip 0.2em}
q_c &=\frac{1}{d}\tr\left[\left(\hat{q}_c\Sigma+ \vec{\hat \mu}\vec{\hat \mu}^\top  \right)\Sigma\hat\Sigma^{-2}\right]\\ \noalign{\vskip 0.2em}
m &= \frac{1}{d}\tr \left[ \vec{\hat \mu}\vec{\theta}_0^\top \Sigma \hat\Sigma^{-1}\right]\\ \noalign{\vskip 0.2em}
h_c &=\frac{1}{d}\tr\left[\vec{\hat \mu}\vec{\mu}_c^\top\hat\Sigma^{-1}
\right] 
\end{cases} 
\end{align}
In particular, these assumptions imply that $V_c, \hat V_c, m_c, \hat m_c$ do not depend on the class labels $c$, and hence
\[ \hat \Sigma = \hat V \Sigma. \]
Now, consider the assignment
\begin{equation}
    \hat{h}_c = 0 \quad \text{and} \quad h_c  = \pi_c;
\end{equation}
is is easy to check that they are a fixed point of the above saddle-point equations, since then
\[ h_c = \frac{\hat m \pi_c}{\hat V} = \pi_c \]
Plugging in this relation in the expression of the update for $\{\hat{q}_c\}_{c \in \mathcal{C}}$ we obtain:
\begin{align}
    \hat{q}_c^\star = \frac{\alpha p_c}{(1+V_c^\star)^2}(\rho_c + \Delta + q_c^\star -2m_c^\star +(h_c^\star-\pi_c)^2) = \frac{\alpha p_c}{(1+V_c^\star)^2}(\rho_c + \Delta + q_c^\star -2m_c^\star)
\end{align} 
Hence by looking at the expression of generalization and training error for ridge regression in eq.~\eqref{eq:app:metrics_definition} we conclude that they will not depend on the values of the mean overlaps $\{h_c^\star\}_{c \in \mathcal{C}}$ and we can prove Gaussian universality in the same way as we did  in the proof of Theorem.~\ref{prop:cov_universality}.
\begin{figure*}
\centering
\includegraphics[width=0.8\textwidth]{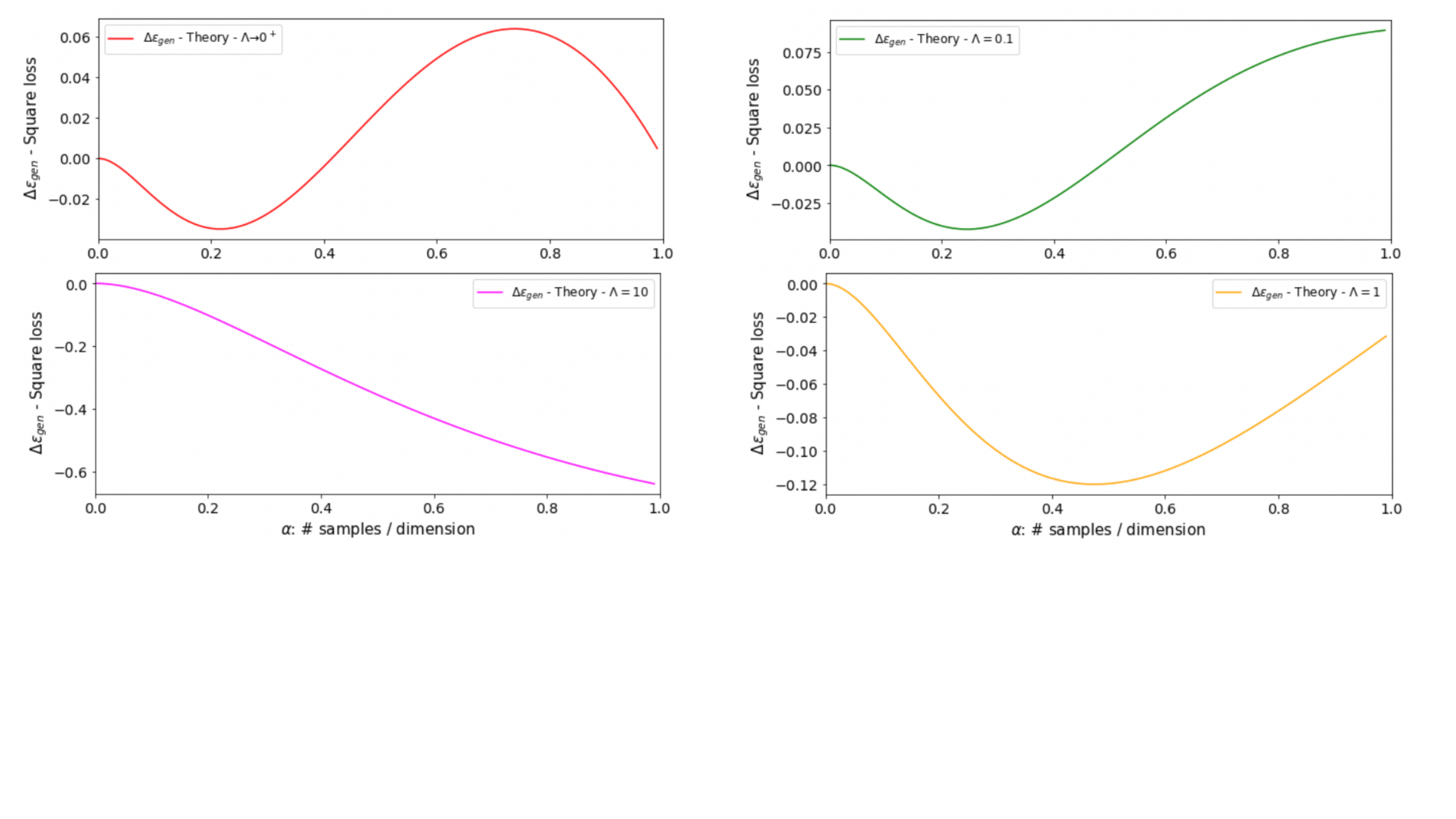}
\vspace{-7em}
\caption{Theoretical prediction from eq.~\eqref{eq:app:theor_discrepanc} for the difference between the generalization error computed from the Gaussian theory and GMM one plotted as a function of the number of samples used in the fit in Algorithm.~\ref{alg:interp_teach}. We use different regularization strength $\Lambda \in \{1e-10,0.1,1,10\}$, increasing clockwise in the figure.}
\label{fig:why_gauss_works}
\end{figure*}
\subsection{Interpolating teachers}
  Different works observed over different datasets that if we find an interpolating teacher vector, in such a way that we can keep the real labels to study ERM performance, the theoretical predictions coming from Gaussian asymptotics would agree with simulations, \cite{gcm_paper} among others. This observation  apparently contradicts the result of Theorem.~\ref{prop:non_univ_mixt} in which we show that correlated labels (as we expect real labels to be) break Gaussian universality. We try to motivate this in the same setting as Theorem.~\ref{prop:non_univ_mixt}:
\begin{theorem} (Real labels universality)
\label{prop:real_labels}
Consider the same setting of Theorem.~\ref{prop:non_univ_mixt}. Fix now the teacher vector to be the maximally correlated one $\Vec{\theta} = \Vec{\mu}$ and find an interpolating teacher $\Tilde{\Vec{\theta}}$ using $\Tilde{n}$ samples coming from the GMM. If we compare the exact asymptotics of Gaussian and GMM for regularization $\Lambda$ at fixed $\alpha = \frac{\Tilde{n}}{d}$ we obtain a discrepancy:
\begin{align}
   \Delta \varepsilon_{\text{gen}} = \frac{4 (E-1) \alpha ^2 \left(E^2 \gamma +2 e-1\right)}{(\alpha +1)^2 (E \gamma +1)^2}
\end{align}
with E($\alpha$,$\Lambda$) solution of:
\begin{align}
    E = 1 + \frac12 \left( \alpha - 1 + \Lambda - \sqrt{4\Lambda + (\alpha -1 + \Lambda)^2} \right). 
\end{align}
\end{theorem}
The result above solve the apparent contradiction: if we strongly overparametrize the fitting model, interpolating teachers can be uncorrelated  with the data structure. 
\begin{algorithm}
\label{alg:interp_teach}
\caption{Get theoretical learning curves keeping real labels}
\begin{algorithmic}
    \STATE  {\bfseries Data generation:} 
    Generate the data from the 2-clusters GMM described in \ref{prop:real_labels} and set:
    \begin{align}
       \Vec{\mu}_+ = - \Vec{\mu}_- = \sqrt{\gamma}(1,\dots,1)^\top, \qquad \Sigma_+ = \Sigma_- = \mat{I}_d, \qquad \Vec{\theta}_0 = \Vec{\mu}_+
    \end{align} 
    \STATE {\bfseries Fit interpolating teacher} Perform ridge regression  for $\lambda \to 0$. Take $\alpha = \frac{n}{d} < 1 $ so that the ERM estimator will interpolate the data.
    \STATE {\bfseries Run exact asymptotics} Compute the theoretical performance of the associated 2-GMM and equivalent GCM from the estimated $\hat{\theta}$ for the same $\alpha$ with a fixed regularization parameter $\Lambda$.
\end{algorithmic}
\end{algorithm}

The proof of the results relies strongly on the fact that we can analyze analytically the performance in this scenario  being it the same as the previous section.  Reminding from eq.~\eqref{eq:no_lamb_eta} that in the limit $\lambda \to 0$ we have in the overparametrized setting $\eta(\alpha,\lambda=0) = \alpha$, we can write: 
\begin{align}
    h_{\textit{erm}} &= \frac{\Vec{\mu}^\top\hat{\Vec{\theta}}}{d} = \gamma \alpha\frac{\gamma + 1}{1+\gamma \alpha} \\ 
    q_{\textit{erm}} &= \frac{\hat{\Vec{\theta}}^\top \Sigma \hat{\Vec{\theta}}}{d} = 
 \alpha  \left(-\frac{\Delta }{\alpha -1}-\frac{(\alpha -1) \pi^2}{\alpha  \gamma +1}+\rho \right) 
\end{align}
We summarize the algorithmic procedure in in Algorithm.~\ref{alg:interp_teach}.
In order to compute the generalization error from eq.~\eqref{eq:gen_error_full generality} we need the quantities $\{\rho,\pi\}$. 
Note that in Algorithm.~\ref{alg:interp_teach} we compute the theoretical curve for the GMM and the GCM from the estimated $\hat{\Vec{\theta}}$. This translates into the fact that the estimated overlaps from numerical simulations in the previous step become the equivalent to: 
\begin{align}
\label{eq:app:new_overlaps}
    \pi = h_{erm} && \rho = q_{erm} 
\end{align}
We fix the regularization for the theoretical prediction to be $\Lambda$, and for simplicity we use the same sample complexity parameter $\alpha$. 
We are now in a position to compare Gaussian and GMM theoretical prediction, more precisely we analyze the difference of the generalization errors for the two data model. Recalling the expression for the generalization error in eq.~\eqref{eq:gen_error_full generality}, we just need to plug the estimated overalps in eq.~\eqref{eq:app:new_overlaps} to obtain:
\begin{align}
\label{eq:app:theor_discrepanc}
   \Delta \varepsilon_{\text{gen}} = \frac{4 \alpha^2 (E-1) \left(E^2 \gamma +2 e-1\right)}{(\alpha +1)^2 (E \gamma +1)^2}
\end{align}
\noindent where $E(\alpha,\Lambda)$ is the solution of eq.~\eqref{eq:eta_expression}, retrieving the result in Theorem.~\ref{prop:real_labels}. 
For the sake of clarity we stated the theorem in a simple setting, indeed once we estimated $\hat{\Vec{\theta}}$ we could have decided to change the value of the sample complexity, as we are now interested in running theory curves. However, to simplify the equations we assumed to run the theoretical prediction for the same value of $\alpha$ as in the ERM fit. We present the results in Fig.~\ref{fig:why_gauss_works} for $\gamma=1$: although the true labels are maximally correlated with the data structure, we see that the Gaussian theoretical predictions with the interpolating teacher $\Vec{\hat{\theta}}$ are very close to the GMM one. We remark that the upper-left plot in Fig.~\ref{fig:fig1} is a nice characterization of Theorem.~\eqref{prop:corr_teach_univers}: as we reach the underparametrized regime we restore Gaussian universality for $\Lambda \to 0^+$ even for correlated teachers.

\section{Details on real dataset simulations}
\label{sec:appendix:numerics}
In this section we report the procedure to create the random regression task on real data, see Algorithm.~\ref{alg:learning_curves}. The implementation of the different numerical simulations described in this work are available in a \href{https://github.com/lucpoisson/GaussianMixtureUniversality}{GitHub respository}.
\vspace{-0.em}
\begin{algorithm}
\label{alg:learning_curves}
\caption{Random teacher regression on real data}
\begin{algorithmic}
    \STATE  {\bfseries Data processing:} 
    Load data and perform the preprocessing step with a transform matrix $F \in \mathbb{R}^{d \times d^\prime}$, with $d^\prime$ dimension of the images in the chosen dataset. We choose for all the figures in this manuscript $d = 2000$. 
    \STATE {\bfseries Match covariance} Compute $\hat{\Sigma} = \frac1n \mat{X}\mat{X}^\top$
    \STATE {\bfseries Create new labels} Forget the real labels associated with the dataset and create new label according to:
    \begin{align}
        y_{\nu} &= \begin{cases}
            &\frac{\Vec{\theta}_0^\top\Vec{x}_{\nu}}{\sqrt{d}} + \sqrt{\Delta} \xi \qquad \text{Ridge regression} \\
            &\rm{sign}\left(\frac{\Vec{\theta}_0^\top\Vec{x}_{\nu} }{\sqrt{d}}+ \sqrt{\Delta} \xi \right)\qquad \text{Logistic regression}
        \end{cases} \\
        \Vec{\theta}_0 &\sim \mathcal{N}(\Vec{0},\mat{I}_d) \qquad \xi \sim \mathcal{N}(0,1)
    \end{align}
    \STATE  {\bfseries Run learning curves} Fix the regularization parameter $\lambda$. 
    \FOR{$\alpha$ in a given range}
    \STATE {\bfseries Simulation} Solve the ERM in eq.~\eqref{eq:main:erm} using sklearn package \textit{LogisticRegression} \cite{scikit-learn} or the Moore-Penrose pseudo-inverse for the ridge estimator:
    \begin{align}
        \hat{\Vec{\theta}} &= \begin{cases}
        &\mat{X}^\top (\mat{X}\mat{X}^\top + \lambda \mat{I}_d)^{-1}Y \qquad n<d\\
        &(\mat{X}^\top\mat{X} + \lambda \mat{I}_d)^{-1}\mat{X}^\top y \qquad \qquad n>d
        \end{cases}
    \end{align}
    \STATE  {\bfseries Gaussian theory} Solve saddle point equations in eqs.~\eqref{eq:saddle:gcm_l2_hat},\eqref{eq:saddle:gcm_l2_nohat} defining Gaussian asymptotics for the model.~\ref{model:gcm}.
    \STATE  {\bfseries Compute erros} Compute the training loss and generalization error using the metrics: 
    \begin{align}
        g(y,\hat{y}) &= \begin{cases}
            (y-\hat{y})^2 \qquad \text{Ridge regression} \\
            \mathbb{P}(y \neq \hat{y}) \qquad \text{Logistic regression}
        \end{cases}
    \end{align}
    \ENDFOR
\end{algorithmic}
\end{algorithm}

\newpage


\end{document}